\documentclass[10pt,fleqn]{article}
\usepackage{latexsym,amsfonts,amssymb,amsbsy,amsmath}
\usepackage{enumerate}
\usepackage{graphicx}

\def\nastavak{eps}

  \newtheorem {lemma}      {Lemma}  [section]

  \newtheorem {theorem}    [lemma]  {Theorem}
  \newtheorem {remark}     [lemma]  {Remark}
  \newtheorem {corollary}  [lemma]  {Corollary}
  \newtheorem {proposition}[lemma]  {Proposition}
  
  \newtheorem {conjecture} [lemma]  {Conjecture}

\newcommand {\proof}[1][{}] {\hbox{\bf Proof #1.\ }}
\newcommand {\qed} {\null \hfill \rule{2mm}{2mm}}

\def\lover#1#2{\lower #1 ex\hbox{$#2$}}

\newcommand{\dogsa}[1]{{\Psi^{#1}_{#1}}}                            %%%
\newcommand{\dogs}[2]{{\Psi^{#1}_{#2}}}                             %%%
\newcommand{\pes}[4]{{\dogs{#1\ldots #2}{#3\ldots #4}}}             %%%
\newcommand{\pesk}[6]{\dogs{#1\ldots\widehat{#2}\ldots #3}{#4\ldots\widehat{#5}\ldots #6}}

\newcommand{\peskd}[5]{{\dogs{#1\ldots #2}{#3\ldots\widehat{#4}\ldots #5}}}
\newcommand{\dogston}{{\dogsa{1\ldots n}}}                          %%%
\newcommand{\dogstonwok}{{\dogsa{1\ldots\widehat{k}\ldots n}}}      %%%
\newcommand{\dogstona}[1]{{\dogsa{1\ldots {#1}}}}                   %%%
\newcommand{\dogstonwoka}[1]{{\dogsa{1\ldots\widehat{k}\ldots {#1}}}} %
\newcommand{\twodots}[1]{{{1}..{#1}}}                               %%%
\def\DJ{{D\kern-.8em\raise.27ex\hbox{-$\!$-}\kern.3em}}

\def\Z{{\mathbb Z}}
\def\Q{{\mathbb Q}}
\def\R{{\mathbb R}}
\def\C{{\mathbb C}}
\def\P{{\mathbb P}}

\title{{\bf Verification and Strengthening of the Atiyah--Sutcliffe Conjectures for Several
Types of Configurations}}

\author{{\bf Dragutin Svrtan}%\thanks{}
\\
\small Department of Mathematics, University of Zagreb,\\[-0.8ex]
\small Bijeni\v{c}ka cesta 30, 10000 Zagreb, Croatia,\\[-0.8ex]
\small \texttt{dsvrtan@math.hr}\\[5mm]
{\bf Igor Urbiha}%\thanks{}
\\
\small Department of Informatics, Polytechnic of Zagreb, University of Zagreb,\\[-0.8ex]
\small Konavoska 2, 10000 Zagreb, Croatia,\\[-0.8ex]
\small \texttt{urbiha@vtszg.hr}}

\date{%\small Submitted: January 1, 2001;  Accepted: January 2, 2001.\\
}

\begin{document}

\maketitle

\begin{abstract}
   In 2001 Sir M. F. Atiyah formulated a conjecture C1 and later with
P. Sutcliffe two stronger conjectures C2 and C3. These conjectures,
inspired by physics (spin-statistics theorem of quantum mechanics),
are geometrically defined for any configuration of points in the
Euclidean three space. The conjecture C1 is proved for $n = 3, 4$
and for general $n$ only for some special configurations (M. F.
Atiyah, M. Eastwood and P. Norbury, D.{\DJ}okovi\'{c}). Interestingly the
conjecture C2 (and also stronger C3) is not yet proven even for
arbitrary four points in a plane. So far we have verified the
conjectures C2 and C3 for parallelograms, cyclic quadrilaterals and
some infinite families of tetrahedra.

   We have also proposed a strengthening of conjecture C3 for configurations
of four points (Four Points Conjectures).

   For almost collinear configurations (with all but one point on a line)
we propose several new conjectures (some for symmetric functions)
which imply C2 and C3. By using computations with multi-Schur
functions we can do verifications up to $n=9$ of our conjectures. We
can also verify stronger conjecture of {\DJ}okovi\' c which imply C2
for his nonplanar configurations with dihedral symmetry.

 \indent Finally we mention that by minimizing  a
geometrically defined energy, figuring in these conjectures, one
gets a connection to some  complicated physical theories, such as
Skyrmions and Fullerenes.
\end{abstract}

%%%%%%%%%%%%%%%%%%%%%%%%%%%%%%%%%%%%%%%%%%%%%%%%%%%%%%%%%%%%%%%%%%%%%%%%%%%%%%%%%%%%%%%%%%%%%%%%%%%%%%%%%%%%

\newcommand{\be}{\begin{equation}}
\newcommand{\ee}{\end{equation}}
\newcommand{\pr}{\partial}
\newcommand{\ie}{{\it i.e.\e }}
\newcommand{\bphi}{\mbox{\boldmath $\phi$}}
\newcommand{\bx}{{\bf x}}
\newcommand{\CP}{\C\P}

\section{Introduction on Geometric Energies}

In this Section we describe some geometric energies,
introduced by Atiyah.
To construct first geometric energy consider $n$ distinct
ordered points, ${\bf x}_i\in\R^3$ for $i=1,...,n$.
 For each pair $i\ne j$ define the unit vector
\be
{\bf v}_{ij}=\frac{\bx_j-\bx_i}{|\bx_j-\bx_i|}
\label{unitv}
\ee
giving the direction of the line joining $\bx_i$ to $\bx_j.$
Now let $t_{ij}\in\CP^1$ be the point on the Riemann sphere associated with
the unit vector ${\bf v}_{ij}$, via the identification
$\CP^1\cong S^2,$ realized as stereographic projection.
Next, set $p_i$ to be the polynomial in $t$ with
roots $t_{ij}$ ($j\ne i$), that is
\be
p_i=\alpha_i\prod_{j\ne i}(t-t_{ij})
\label{defp}
\ee
where $\alpha_i$ is a certain normalization coefficient.
In this way we have constructed $n$ polynomials which all have degree
$n-1,$ and so we may write
\[
p_i=\sum_{j=1}^n m_{ij}t^{j-1}.
\]
Finally, let $M_n$ be the
$n\times n$ matrix with entries $m_{ij},$ and let ${\mathcal D}_n$ be its
determinant
\be
{\mathcal D}_n={\mathcal D}_n({\bf x}_1,...,{\bf x}_n)=\mbox{det}\ M_n.
\ee
This geometrical construction is relevant to the Berry-Robbins
problem, which is concerned with specifying how a spin
basis varies as $n$ point particles move in space, and supplies
a solution provided it can be shown that ${\mathcal D}_n$ is always non-zero.
For $n=2,3,4$ it can be proved that ${\mathcal D}_n\ne 0$ (Atiyah $n=3$, Eastwood and Norbury $n=4$) and
numerical computations suggest that $|{\mathcal D}_n|\ge 1$ for all $n,$
with the minimal value $|{\mathcal D}_n|=1$ being attained by $n$ collinear points.

The geometric energy is the $n$-point energy defined by
\be
E_n=-\log |{\mathcal D}_n|,
\label{geom}
\ee
so minimal energy configurations maximize the modulus of the determinant.

This energy is geometrical in the sense that it only depends on
the directions of the lines joining the points, so it is
translation, rotation and scale invariant. Remarkably, the minimal
energy configurations, studied numerically for all $n\le 32,$ are
essentially the same as those for the Thomson problem.

%%%%%%%%%%%%%%%%%%%%%%%%%%%%%%%%%%%%%%%%%%%%%%%%%%%%%%%%%%%%%%%%%%%%%%%%%%%%%%%%%%%%%%%

%%%%%%%%%%%%%%%%%%%%%%%%%%%%%%%%%%%%%%%%%%%%%%%%%%%%%%%%%%%%%%%%%%%%%%%%%%%%%%%%%%%%%%%%%%%%%%%%%%%%%%%%%% [1]
\section{Eastwood--Norbury formulas for Atiyah determinants}

In this section we first recall Eastwood--Norbury formula for
Atiyah determinant for three or four points in Euclidean
three--space. In the case $n=3$ the (non normalized) Atiyah determinant reads as
\[
D_3=d_3(r_{12},r_{13},r_{23})+8r_{12}r_{13}r_{23}
\]
where \
\[
d_3(a,b,c)=(a+b-c)(b+c-a)(c+a-b)
\]
and $r_{ij}$ $(1\leq i<j\leq 3)$ is the distance
between the $i^{\mbox{th}}$ and $j^{\mbox{th}}$ point.

The normalized Atiyah determinant for $3$ points is
\[
{\mathcal D}_3=\frac{D_3}{8r_{12}r_{13}r_{23}}
\]
and it is evident that $|{\mathcal D}_3|={\mathcal D}_3\geq 1$.

In the case $n=4$ the (non normalized) Atiyah determinant $D_4$ has real part
given by a polynomial (with $248$ terms) as follows:
\begin{equation}
\label{E:1.0}
\Re (D_4)=64r_{12}r_{13}r_{23}r_{14}r_{24}r_{34}-4d_3(r_{12}r_{34},r_{13}r_{24},r_{14}r_{23})+A_4+288V^2
\end{equation}
where
\[
A_4=\sum_{l=1}^{4}\left(\sum_{\scriptstyle (l\neq) i=1}^{4}
r_{li}((r_{lj}+r_{lk})^2-r_{jk}^2)\right)d_3(r_{ij},r_{ik},r_{jk})
\]
(here $\{j,k\}=\{1,2,3,4\}\setminus\{l,i\}$) and $V$ denotes the
volume of the tetrahedron with vertices our four points:
\begin{equation}
\label{E:144V}
\begin{array}{ll}
144V^2 =
&
r_{12}^2r_{34}^2(r_{13}^2+r_{14}^2+r_{23}^2+r_{24}^2-r_{12}^2-r_{34}^2)+\
\lover{-0.2}{\mbox{{\tiny two similar terms}}}\\[1ex]
&-(r_{12}^2r_{13}^2r_{23}^2+\lover{-0.2}{\mbox{{\tiny three similar terms}}})\\[1ex]
\end{array}
\end{equation}
We now state two formulas which will be used later:
\begin{enumerate}
    \item Alternative form of $A_4$:
\begin{equation}
\label{E:A4-a}
\begin{array}{ll}
A_4=
&\displaystyle
 \sum_{ {\scriptstyle l=1 } }^4
\left(( d_3(r_{il},r_{jl},r_{kl})+8r_{il}r_{jl}r_{kl}
+r_{il}(r_{il}^2-r_{jk}^2) +\right.\\
&
\hphantom{ \sum_{ {\scriptstyle l=1 } }^4}
r_{jl}(r_{jl}^2-r_{ik}^2)+\left. r_{kl}(r_{kl}^2-r_{ij}^2)\right)
d_3(r_{ij},r_{ik},r_{jk}),
\end{array}
\end{equation}
where for each $l$ we write $\{1,2,3,4\}\setminus\{l\}=\{i<j<k\}$.
    \item The sum of the second and the fourth term of (\ref{E:1.0}) can be rewritten as
\begin{equation}
\label{E:e0}
\begin{array}{l}
144V^2-2d_3(r_{12}r_{34},r_{13}r_{24},r_{14}r_{23})=\\
\hspace{0.7cm}\begin{array}{l}
=(r_{12}-r_{34})^2(r_{13}^2r_{24}^2+r_{14}^2r_{23}^2-r_{12}^2r_{34}^2)+\
\lover{-0.2}{\mbox{{\tiny two such terms}}}\ +\\
\hphantom{=}+4r_{12}r_{13}r_{23}r_{14}r_{24}r_{34}-\\
\hphantom{=}-r_{12}^2r_{13}^2r_{23}^2-r_{12}^2r_{14}^2r_{24}^2
-r_{13}^2r_{14}^2r_{34}^2-r_{23}^2r_{24}^2r_{34}^2.
\end{array}
\end{array}
\end{equation}
\end{enumerate}
It is well known that this quantity is always nonpositive.

The imaginary part $\Im (D_4))$ of Atiyah determinant can
be written as a product of $144V^2$ with a polynomial (with
integer coefficients) having 369 terms.

The normalized Atiyah determinant for $4$ points is
\[
{\mathcal D}_4=\frac{D_4}{\displaystyle 2^{\binom{4}{2}}\prod_{1\leq i<j\leq 4} r_{ij}}.
\]

The original Atiyah conjecture in our cases is equivalent to nonvanishing
of the determinants $D_3$ and $D_4$.

A stronger conjecture of Atiyah and Sutcliffe
(\cite{AS},Conjecture 2) states in our cases that $|D_3|$ $\geq$
$8r_{12}r_{13}r_{23}$ ($\Leftrightarrow |{\mathcal D}_3|\geq 1$) and $|D_4|\geq
64r_{12}r_{13}r_{23}r_{14}r_{24}r_{34}$ ($\Leftrightarrow |{\mathcal D}_4|\geq 1$).

From the formula (\ref{E:1.0}) above, with the help of the simple inequality
$d_3(a,b,c)$ $\leq$ $abc$ (for $a,b,c\geq 0$), Eastwood and Norbury got
"almost" the proof of the stronger conjecture by exhibiting the inequality
\[
\Re(D_4)\geq 60r_{12}r_{13}r_{23}r_{14}r_{24}r_{34}.
\]
To remove the word "almost" seems to be not so easy (at present
not yet done even for planar configuration of four points).

%%%%%%%%%%%%%%%%%%%%%%%%%%%%%%%%%%%%%%%%%%%%%%%%%%%%%%%%%%%%%%%%%%%%%%%%%%%%%%%%%%%%%%%%%%%%%%%%%%%%%%%%%% [2]
A third conjecture (stronger than the second) of Atiyah and
Sutcliffe (\cite{AS}, Conjecture 3) can be expressed, in the four
point case, in terms of polynomials in the edge lengths as
\begin{equation}
\label{E:e1}
|D_4|^2\geq
\prod_{\{i<j<k\}\subset\{1,2,3,4\}}(d_3(r_{ij},r_{ik},r_{jk})+8r_{ij}r_{ik}r_{jk})
\end{equation}
where the product runs over the four faces of the tetrahedron.\\
(cf.\ {\tt
ftp://ftp.maths.adelaide.edu.au/pure/meastwood/atiyah.ps})

%%%%%%%%%%%%%%%%%%%%%%%%%%%%%%%%%%%%%%%%%%%%%%%%%%%%%%%%
%
%  Added on 20060603
%
%%%%%%%%%%%%%%%%%%%%%%%%%%%%%%%%%%%%%%%%%%%%%%%%%%%%%%%%

In the first part of this paper we study some infinite families of quadrilaterals and tetrahedra and
verify both Atiyah and Sutcliffe conjectures for several such infinite families.
In this version of the paper we propose a somewhat stronger conjecture than (\ref{E:e1})
which reads as follows:\\

\begin{conjecture}(Four Points Conjectures)\\
\begin{equation}
\label{OC:1}
\begin{array}{l}
\displaystyle \Re(D_4)-(4+\mbox{$\frac{3}{4}$})\cdot 288V^2\geq\\[2mm]
\displaystyle \geq 64\prod_{1\leq i,j\leq 4}r_{ij}+
\sum_{\{i<j<k\}\subset\{1,2,3,4\}}(4+\mbox{$\frac{1}{4}$}\delta)r_{il}r_{jl}r_{kl}
d_3(r_{ij},r_{ik},r_{jk})
\end{array}
\end{equation}
where
\[
\delta=\left\{
\begin{array}{c@{\ ,\ }l}
\displaystyle
\frac{d_3(r_{ij},r_{ik},r_{jk})}{r_{ij}r_{ik}r_{jk}} & \mbox{weak version}\\[4mm]
1 &  \mbox{strong version}
\end{array}
\right.
\]
\end{conjecture}
\begin{proposition}
Any of the Four Points Conjectures (\ref{OC:1}) imply conjecture (\ref{E:e1}).
\end{proposition}
\proof By using the inequality $1\geq d_3(a,b,c)/(abc), (a,b,c>0)$ (see Appendix 2,
Proposition \ref{P:A2}) we see that the strong version implies the weak version of conjecture.
We then rewrite the rhs of the weak version of (\ref{OC:1}) as follows:
\[
\prod_{1\leq i,j\leq 4}r_{ij}
\left(
\frac{1}{4}
\sum_{l=1}^4
\left(
\frac{d_3(r_{ij},r_{ik},r_{jk})}{r_{ij}r_{ik}r_{jk}}+8
\right)^2
\right)
\]
Finally, by the quadratic--geometric (QG) inequality we obtain
\[
\geq
\prod_{1\leq i,j\leq 4}r_{ij}
\left(
\prod_{l=1}^4
\left(
\frac{d_3(r_{ij},r_{ik},r_{jk})}{r_{ij}r_{ik}r_{jk}}+8
\right)
\right)^\frac{2}{4}
=
\left(
\prod_{l=1}^4
\left(
d_3(r_{ij},r_{ik},r_{jk})+8r_{ij}r_{ik}r_{jk}
\right)
\right)^\frac{1}{2}
\]
Thus we obtain:
\[
|D_4|^2\geq
|\Re (D_4)|^2\geq
|\Re (D_4)-(4+\mbox{$\frac{3}{4}$})\cdot 288V^2|^2\geq
\prod_{l=1}^4
(d_3(r_{ij}r_{ik}r_{jk})+8r_{ij}r_{ik}r_{jk})
\]
i.e.\ the inequality (\ref{E:e1}).\qed

\begin{remark}
In terms of trigonometry
(see subsection "Atiyah determinant for triangles and quadrilaterals via trigonometry"
on page \pageref{S:Adftaqvt}),
the weak Four Points Conjecture can be written simply as
\[
\Re(D_4)-(4+\mbox{$\frac{3}{4}$})\cdot 288V^2
\geq
\left(
\prod_{1\leq i,j\leq 4}r_{ij}
\right)
\left(
4
\sum_{l=1}^4
c_l^2
\right)
\]
where
\[
c_l:=\cos^2\frac{X^{(l)}}{2}+\cos^2\frac{Y^{(l)}}{2}+\cos^2\frac{Z^{(l)}}{2},\ \ l=1,2,3,4.
\]
and $X^{(l)}$, $Y^{(l)}$, $Z^{(l)}$ are the angles of the triangle opposite to the vertex $l$.
\end{remark}

%%%%%%%%%%%%%%%%%%%%%%%%%%%%%%%%%%%%%%%%%%%%%%%%%%%%%%%%
%
%  End of Added on 20060603
%
%%%%%%%%%%%%%%%%%%%%%%%%%%%%%%%%%%%%%%%%%%%%%%%%%%%%%%%%

\subsection{Atiyah--Sutcliffe conjecture for (vertically) upright tetrahedra (or pyramids)}

We call a tetrahedron upright if some of its vertices (say 4) is
equidistant from all the remaining vertices (1, 2 and 3, which we
can think as lying in a horizontal plane.)
\begin{center}
\begin{tabular}{cl}
\begin{minipage}{125pt}
\centerline{\includegraphics[width=125pt,keepaspectratio]{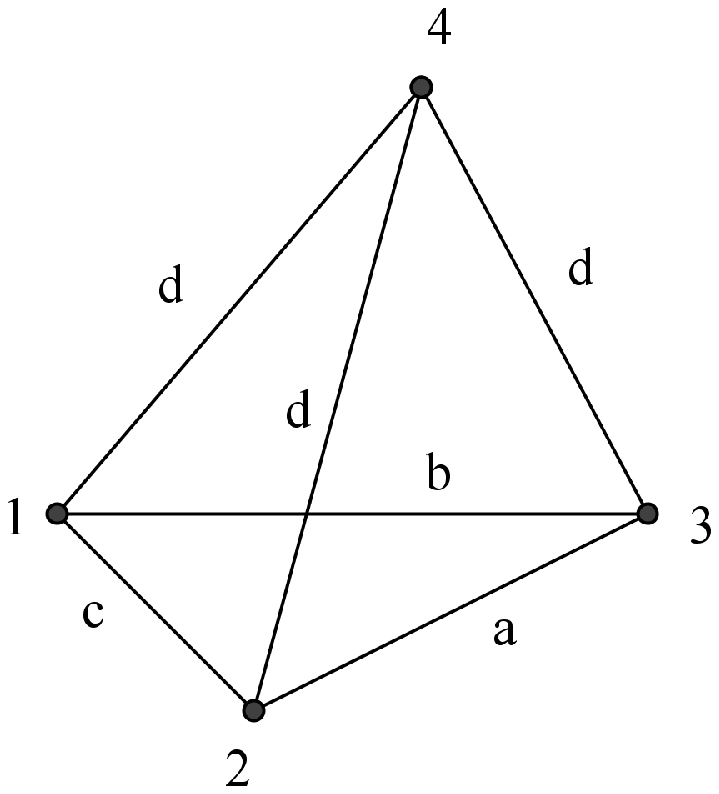}}
\end{minipage}
&
\begin{minipage}{200pt}
$r_{23}=a$\\
$r_{13}=b$\\
$r_{12}=c$\\
$r_{14}=r_{24}=r_{34}=d$
\end{minipage}
\end{tabular}
\end{center}

Note that then $d\geq R=$ the circumradius of the base triangle 123, then
by Heron's formula we have: $R={abc}/{\sqrt{(a+b+c)d_3(a,b,c)}}$.

Here, as before, $d_3(a,b,c)=(a+b-c)(a-b+c)(-a+b+c)$, $(a,b,c>0)$.

The left hand side of the strong Four Points Conjecture \ref{OC:1} (but without $\frac{3}{4}$ term!)
can be evaluated as follows, by using Eastwood-Norbury formula (\ref{E:1.0})
\[
\begin{array}{l@{\ =\ }l}
LHS&\Re(D_4)-4\cdot 288V^2=\\
&\displaystyle
64\!\!\!\!\!\prod_{1\leq i<j\leq 4} r_{ij}-4d_3(r_{12}r_{34},r_{13}r_{24},r_{14}r_{23})
+A_4-3\cdot 288V^2
\end{array}
\]
where
\[
-4d_3(r_{12}r_{34},r_{13}r_{24},r_{14}r_{23})=-4d_3(a,b,c)d^3
\]
\[
\begin{array}{ll}
A_4 = & \displaystyle
\sum_{l=1}^4\left(\sum_{(l\neq)i=1}^4r_{li}((r_{lj}+r_{lk})^2-r_{jk}^2)d_3(r_{ij},r_{ik},r_{jk})\right)=\\[5mm]
\hphantom{A_4 }= & \displaystyle
\!\!\!\!\!\!\sum_{cyc(a,b,c)}\left[c((b+d)^2-d^2)+b((c+d)^2-d^2)+d((b+c)^2-a^2)\right]d_3(a,d,d)+\\[5mm]
&+ \left[d((d+d)^2-a^2)+d((d+d)^2-b^2)+d((d+d)^2-c^2)\right]d_3(a,b,c)=\\[5mm]
\hphantom{A_4 }= & \displaystyle
[4bcd+((b+c)^2-a^2)d+b^2c+bc^2](2a^2d-a^3)]+\cdots+\\[5mm]\overline{}
& \displaystyle
+[12d^3-(a^2+b^2+c^2)d]d_3(a,b,c)\hspace{2cm}(by\ \  \ref{E:A4-a})
\end{array}
\]
\[
\begin{array}{l}
-3\cdot 288V^2=\\
=-6[(b^2+c^2-a^2)a^2d^2+(c^2+a^2-b^2)b^2d^2+(a^2+b^2-c^2)c^2d^2-a^2b^2c^2] (by\ \ref{E:144V})
\end{array}
\]
Similarly the right hand side of the Conjecture \ref{OC:1}
\[
\begin{array}{rl}
RHS =&\displaystyle 64\prod_{1\leq i<j\leq 4}r_{ij}+\sum_{l=1}^4\left(4+\frac{1}{4}\right)r_{il}r_{jl}r_{kl}d_3(r_{ij},r_{ik},r_{jk})\\[5mm]
=&\displaystyle 64abcd^3 +\left(4+\frac{1}{4}\right)bcd(2a^2d-a^3)+\mbox{ two such terms}+\\[5mm]
&+\displaystyle \left(4+\frac{1}{4}\right)d^3d_3(a,b,c)
\end{array}
\]
Now we can rewrite the difference
\[
LHS-RHS=I+II
\]
where
\[
I=\sum_{cyc}(b^2c+bc^2)(2a^2d-a^3)-(a^2+b^2+c^2)d_3(a,b,c)d+6a^2b^2c^2-24\frac{a^2b^2c^2}{a+b+c}d
\]
and
\[
\begin{array}{ll}
II=&\displaystyle
%-\frac{3}{2}[(b^2+c^2-a^2)a^2d^2+(c^2+a^2-b^2)b^2d^2+(a^2+b^2-c^2)c^2d^2-a^2b^2c^2]+\\
%&\displaystyle
\left(4-\frac{1}{4}\right)d_3(a,b,c)d^3+\sum_{cyc}((b+c)^2-a^2)d(2a^2d-a^3)-\\[5mm]
&\displaystyle
-6\sum_{cyc}(b^2+c^2-a^2)a^2d^2-\frac{1}{4}\sum_{cyc}bcd(2a^2d-a^3)+24\frac{a^2b^2c^2}{a+b+c}d
\end{array}
\]
Then we can further simplify
\[
\begin{array}{ll}
I=&\displaystyle
\left[
4abc(ab+ac+bc)-(a^2+b^2+c^2)d_3(a,b,c)-\frac{24a^2b^2c^2}{a+b+c}
\right]d+\\
&\displaystyle
+6a^2b^2c^2-\sum_{sym}a^3b^2c
\end{array}
\]
and
\[
\begin{array}{l}
\displaystyle
II=
d\left[\frac{15}{4}d_3(a,b,c)d^2+(a+b+c)
\left(
\frac{7}{2}abc-4d_3(a,b,c)
\right)d+
24\frac{a^2b^2c^2}{a+b+c}
+
\right.\\
\displaystyle
\hphantom{II=}
\left.
+
\frac{1}{4}
abc(a^2+b^2+c^2)-
(a+b+c)
\left(
\sum_{sym}a^3b-a^4-b^4-c^4
\right)
\right]
\end{array}
\]
\begin{lemma}
\label{L:1a}
We have the following strengthening of the basic inequality for our function
$d_3(a,b,c)=(a+b-c)(a-b+c)(-a+b+c)$:
\[
d_3(a,b,c)\leq \frac{9a^2b^2c^2}{(a+b+c)(a^2+b^2+c^2)}\ \ (\leq \frac{27a^2b^2c^2}{(a+b+c)^3}\leq abc)
\]
\end{lemma}
\proof
We have
\[
\begin{array}{l}
9a^2b^2c^2-(a^2+b^2+c^2)(a+b+c)d_3(a,b,c)=\\
=9a^2b^2c^2-(a^2+b^2+c^2)(2a^2b^2+2a^2c^2+2b^2c^2-a^4+b^4+c^4)=\\
=3a^2b^2c^2-a^4b^2-a^2b^4-a^4c^2-a^2c^4-b^4c^2-b^2c^4+a^6+b^6+c^6=\\
(a^2-b^2)[a^2(a^2-c^2)-b^2(b^2-c^2)]+c^2(a^2-c^2)(b^2-c^2)\geq 0\\
(\mbox{if we assume }a\geq b\geq c\geq 0)
\end{array}
\]
(a special instance of a Schur inequality)\qed\\
(Note that this result follows from the formula $OG^2=R^2-(a^2+b^2+c^2)/9$ for the
distance of the circumcenter and the centroid of a triangle.)

Now we have
\begin{lemma}
\label{L:2a}
The quantity $I$ is increasing w.r.t.\ $d$ and it is positive for $d\geq R$.
\end{lemma}
\proof
We prove that the coefficient of $d$ in $I$ is positive by using that
$(ab+ac+bc)(a+b+c)\geq 9abc$ and  Lemma \ref{L:1a} .

The proof of positivity of $I$ reduces to the positivity of the following quantity:
\[
\begin{array}{l}
\{[4abc(ab+ac+bc)-(a^2+b^2+c^2)d_3(a,b,c)](a+b+c)-24a^2b^2c^2\}^2-\\
-(a+b+c)^3(a^2b+ab^2+a^2c+ac^2+b^2c+bc^2-6abc)^2d_3(a,b,c)
\end{array}
\]
which by substituting $a=b+h$ and $b=c+k$ and then expanding has all coefficients positive
(and ranging from $1$  to $32151$).\qed

\begin{lemma}
\label{L:2b}
The quantity $II$ is increasing w.r.t.\ $d$ and it is positive for $d\geq R$.
\end{lemma}
\proof
Let $II=d\cdot III$. Then
\[
\begin{array}{l}
\frac{\partial III}{\partial d}=\left(\frac{15}{2}d_3(a,b,c)d-\frac{a+b+c}{2}d_3(a,b,c)\right)
+\frac{7(a+b+c)}{2}(abc-d_3(a,b,c))%\\
\end{array}
\]
The second term is positive by Proposition \ref{P:A2}. For the first term we have:
\[
\begin{array}{l}
\frac{15}{2}d_3(a,b,c)d-\frac{a+b+c}{2}d_3(a,b,c)\geq
\frac{15}{2}d_3(a,b,c)R-\frac{a+b+c}{2}d_3(a,b,c)\geq\\[3mm]
(\frac{15abc}{(a+b+c)^{3/2}}-\sqrt{d_3(a,b,c)})\frac{a+b+c}{2}\sqrt{d_3(a,b,c)}\geq 0
\end{array}
\]
by Lemma (\ref{L:1a}).

\noindent
The proof of positivity of $II$ reduces to the positivity of the following quantity:
\[
\begin{array}{l}
\displaystyle
\frac{15}{4}d_3(a,b,c)R^2+(a+b+c)
\left(
\frac{7}{2}abc-4d_3(a,b,c)
\right)R+
24\frac{a^2b^2c^2}{a+b+c}
+\\
+
\frac{1}{4}
abc(a^2+b^2+c^2)-
(a+b+c)
\left(
\sum_{sym}a^3b-a^4-b^4-c^4
\right)
\end{array}
\]
which can be nicely visualized by {\tt Maple} using tangential coordinates ($a=v+w$, $b=u+w$, $c=u+v$).
\qed

%%%%%%%%%%%%%%%%%%%%%%%%%%%%%%%%%%%%%%%%%%%%%%%%%%%%%%%%%%%%%%%%%%%%%%%%%%%%%%%
\subsection{Atiyah--Sutcliffe conjectures for edge--tangential tetrahedra}

By edge--tangential tetrahedron we shall mean any tetrahedron for which there exists a sphere touching
all its edges (i.e.\ its 1--skeleton has an inscribed sphere.) For each $i$ from $1$ to $4$ we denote by
$t_i$ the length of the segment (lying on the tangent line) with one endpoint the vertex and the other
the point of contact of the tangent line with a sphere.

\begin{center}
\begin{tabular}{cl}
\begin{minipage}{100pt}
\centerline{\includegraphics[width=150pt,keepaspectratio]{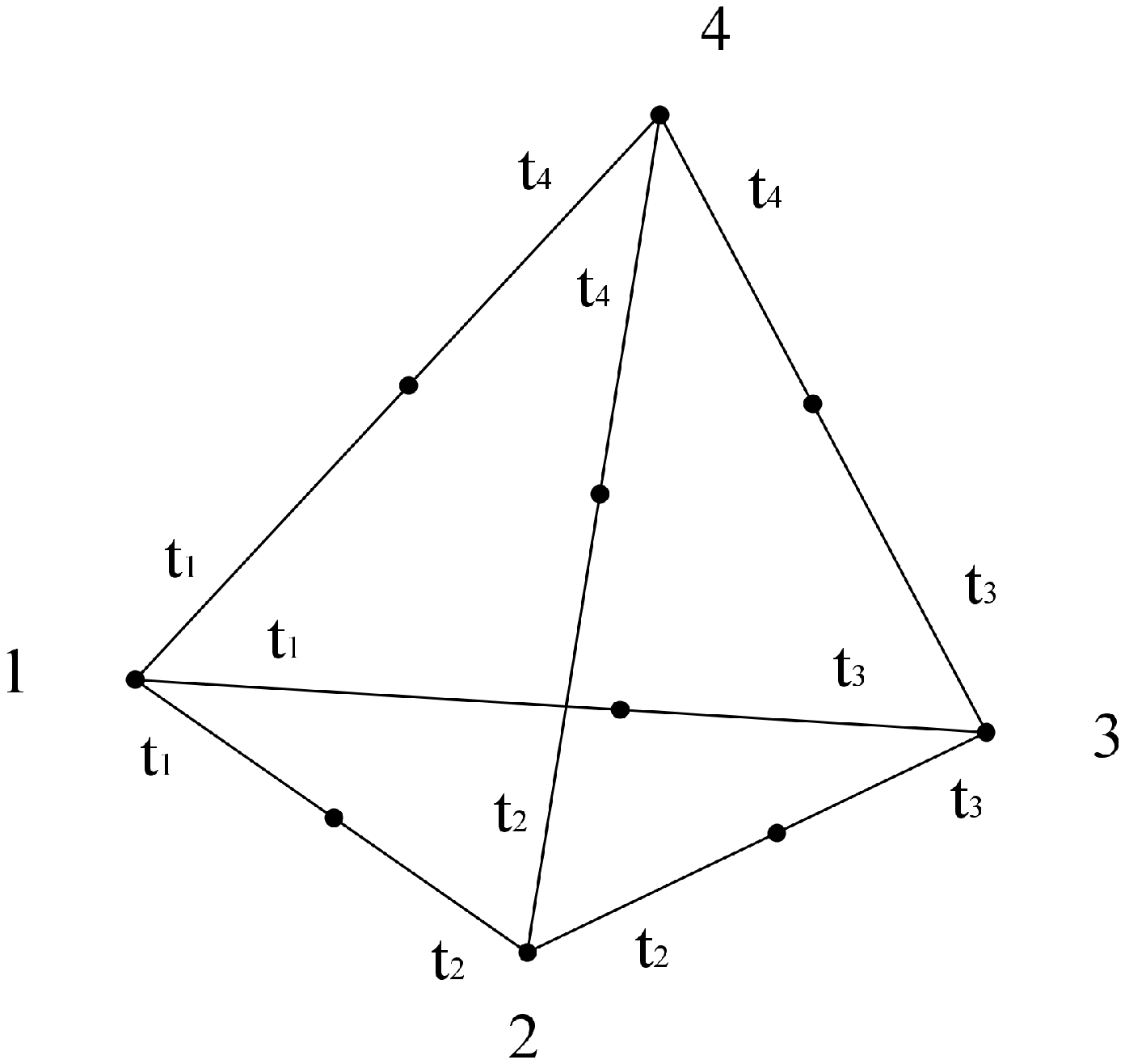}}
\end{minipage}
&
\begin{minipage}{200pt}
\begin{center}
$r_{ij}=t_i+t_j,\ (1\leq i<j\leq 4)$
\end{center}
\end{minipage}
\end{tabular}
\end{center}
Now we shall compute all the ingredients appearing in the Eastwood--Norbury formula for $D_4$
%(the real part of the Atiyah determinant)
in terms of elementary symmetric functions of the (tangential) variables $t_1, t_2, t_3, t_4$
(recall $e_1=t_1+t_2+t_3+t_4$, $e_2=t_1t_2+t_1t_3+t_1t_4+t_2t_3+t_2t_4+t_3t_4$,
$e_3=t_1t_2t_3+t_1t_2t_4+t_1t_3t_4+t_2t_3t_4$, $e_4=t_1t_2t_3t_4$).
\[
\begin{array}{l}
  64r_{12}r_{13}r_{23}r_{14}r_{24}r_{34} =
  \displaystyle 64\hspace{-2ex}\prod_{1\leq i<j\leq 4}\hspace{-2ex}(t_i+t_j)=
  64s_{3,2,1}=\\[6mm]
   \hphantom{64r_{12}r_{13}r_{23}r_{14}r_{24}r_{34}} =64\left|
\begin{array}{ccc}
  e_3 & e_4 & 0 \\
  e_1 & e_2 & e_3 \\
  0 & 1 & e_1
\end{array}
  \right|=
  64e_3e_2e_1-64e_4e_1^2-64e_3^2
\end{array}
\]
Here we have used Jacobi--Trudi formula for the triangular Schur function $s_{3,2,1}$
(see \cite{Jacobi}, (3.5)). Furthermore we have
\[
\begin{array}{r@{\ =\ }l}
  -4d_3(r_{12}r_{34},r_{13}r_{24},r_{14}r_{23}) & 128e_4e_2-32e_4e_1^2-32e_3^2 \\[2mm]
  288V^2 & 128e_4e_2-32e_3^2
\end{array}
\]
In order to compute $A_4$ we first compute, for fixed $l$ the following quantities
\[
\begin{array}{r@{\ =\ }l}
  d_3(r_{ij},r_{ik},r_{jk})  & 8t_it_jt_k\\
  \displaystyle \sum_{(l\neq) i=1}^4r_{li}((r_{lj}+r_{lk})^2-r_{jk}^2)
  & 4(3t_l(t_1+t_2+t_3+t_4)+2(t_it_j+t_it_k+t_jt_k))t_l.
\end{array}
\]
Thus we get:
\[
A_4=32(3e_1^2+4e_2)e_4=96e_4e_1^2+128e_4e_2.
\]
Now we adjust terms in $D_4$, in order to get shorter expression, as follows
\[
\begin{array}{rl}
  D_4 =& (64r_{12}r_{13}r_{23}r_{14}r_{24}r_{34}-2\cdot 288V^2)+\\
                & +(-4d_3(r_{12}r_{34},r_{13}r_{24},r_{14}r_{23})-288V^2)+A_4+4\cdot 288V^2 \\[2mm]
   =& (64e_3e_2e_1-64e_4e_1^2-256e_4e_2)+(-32e_4e_1^2)+\\
   & +(96e_4e_1^2+128e_4e_2)+ 4\cdot 288V^2\\[2mm]
   =&  64e_3e_2e_1-128e_4e_2+1152V^2\\[2mm]
   =&  64e_2(e_3e_1-2e_4)+1152V^2\\[2mm]
   =&  64e_2(2e_4+m_{211})+1152V^2,
\end{array}
\]
where $m_{211}=t_1^2t_2t_3+\cdots$ denotes the monomial symmetric function associated
to the partition $(2,1,1)$.

In order to verify the third conjecture of Atiyah and Sutcliffe
\[
|D_4|^2\geq \prod_{\{i<j<k\}\subset\{1,2,3,4\}}(d_3(r_{ij},r_{ik},r_{jk})+8r_{ij}r_{ik}r_{jk})
\]
we note first that
\[
\begin{array}{l@{\ =\ }l}
  d_3(r_{ij},r_{ik},r_{jk})+8r_{ij}r_{ik}r_{jk} & (8t_it_jt_k+8(t_i+t_j)(t_i+t_k)(t_j+t_k)) \\[2mm]
   & 8(t_i+t_j+t_k)(t_it_j+t_it_k+t_jt_k)
\end{array}
\]
and state the following:
\begin{lemma}
\label{L:2}
For any nonnegative real numbers $t_1,t_2,t_3,t_4\geq 0$ the following inequality
\begin{equation}
\label{E:e2}
\begin{array}{l}
(t_1t_2+t_1t_3+t_1t_4+t_2t_3+t_2t_4+t_3t_4)^2(2t_1t_2t_3t_4+m_{211}(t_1,t_2,t_3,t_4))^2\geq\\[2mm]
\geq \prod_{\{i<j<k\}\subset\{1,2,3,4\}}(t_i+t_j+t_k)(t_it_j+t_it_k+t_jt_k)
\end{array}
\end{equation}
holds true.
\end{lemma}
\proof[of Lemma \ref{L:2}]
The difference between the left hand side and the right hand side of the above inequality
(\ref{E:e2}), written in terms of monomial symmetric functions is equal to
\[
\begin{array}{ll}
LHS-RHS=&m_{6321}+3m_{6222}+m_{543}+2m_{5421}+7m_{5322}+5m_{5331}+\\[2mm]
&+3m_{444}+7m_{4431}+8m_{4422}+8m_{4332}+3m_{3333}\geq 0
\end{array}
\]
\qed
\begin{remark}
One may think that the inequality in Lemma \ref{L:2} can be obtained as a product of two
simpler inequalities. This is not the case, because the following inequalities hold true:
{\small
\[
\begin{array}{rcl}
(t_1t_2+t_1t_3+t_1t_4+t_2t_3+t_2t_4+t_3t_4)^2
&\leq&
 \displaystyle\prod_{\{i<j<k\}\subset\{1,2,3,4\}}(t_i+t_j+t_k) \\[6mm]
(2t_1t_2t_3t_4+m_{211}(t_1,t_2,t_3,t_4))^2
&\geq&
 \displaystyle \prod_{\{i<j<k\}\subset\{1,2,3,4\}}(t_it_j+t_it_k+t_jt_k)  \\
    &
\end{array}
\]
}
\end{remark}
Now we continue with verification of the third conjecture of Atiyah and Sutcliffe for
edge tangential tetrahedron:
\[
\begin{array}{rl}
  |D_4|^2\geq & (D_4)^2\geq [64e_2(2e_4+m_{211})]^2 \\[2mm]
  \geq & \displaystyle 8^4\hspace{-5ex}\prod_{\{i<j<k\}\subset\{1,2,3,4\}}
  \hspace{-4ex}(t_i+t_j+t_k)(t_it_j+t_it_k+t_jt_k)\ \ \ (\mbox{by Lemma \ref{L:2}}) \\[6mm]
  = & \displaystyle\hspace{-5ex}\prod_{\{i<j<k\}\subset\{1,2,3,4\}}
  \hspace{-4ex}(d_3(r_{ij},r_{ik},r_{jk})+8r_{ij}r_{ik}r_{jk})
\end{array}
\]
so the strongest Atiyah--Sutcliffe conjecture is verified for edge--tangential tetrahedra.

\subsection{Verification of the strong Four Points Conjecture for edge--tangential tetrahedra}

The strong Four Points Conjecture \ref{OC:1} for edge tangential tetrahedra is equivalent to positivity
of the following quantity:
\[
\begin{array}{l}
\displaystyle
\Re(D_4)-64\prod r_{ij}-(4+\mbox{$\frac{3}{4}$})288V^2-
\sum_{l=1}^4(4+\mbox{$\frac{1}{4}$})r_{il}r_{jl}r_{kl}\, d_3(r_{ij},r_{ik},r_{jk})\\
=(-d_3(r_{12}r_{34},r_{13}r_{24},r_{14}r_{23})+A_4)+288V^2-(4+\mbox{$\frac{3}{4}$})288V^2\\
\displaystyle
\hphantom{=}
-\sum_{l=1}^4(4+\mbox{$\frac{1}{4}$})r_{il}r_{jl}r_{kl}\, d_3(r_{ij},r_{ik},r_{jk})\\
\displaystyle
=(-32m_{3111}-32m_{222}+96m_{3111}+320m_{2211})-240m_{2211}+120m_{222}\\
\displaystyle
\hphantom{=}
-(34m_{3111}+136m_{2211}+34m_{222})\\
\displaystyle
=
30m_{3111}+54m_{222}-56m_{2211}
\end{array}
\]
In terms of augmented monomial symmetric functions
\[
\widetilde{m}_\lambda (t_1,t_2,t_3,t_4)= \sum_{\sigma\in S_4}t^{\sigma.\lambda}
\]
the last quantity is equal to
\[
=5\widetilde{m}_{3111}+9\widetilde{m}_{222}-14\widetilde{m}_{2211}
\ (\geq 0\mbox{ by Muirheads's inequality})
\]
Thus, the strong Four Points Conjecture is verified for the edge--tangential tetrahedra.

Note that the verification of this conjecture which is stronger than A--S conjecture C3 is
somewhat simpler (at least for edge--tangential tetrahedra).
%%%%%%%%%%%%%%%%%%%%%%%%%%%%%%%%%%%%%%%%%%%%%%%%%%%%%%%%%%%%%%%%%%%%%%%%%%%%%%
%
%     20060803
%
%%%%%%%%%%%%%%%%%%%%%%%%%%%%%%%%%%%%%%%%%%%%%%%%%%%%%%%%%%%%%%%%%%%%%%%%%%%%%%

\subsection{Trirectangular tetrahedra}

A tetrahedron is called trirectangular if it has a vertex at which all the face
angles are right angles. The opposite face to such a vertex we call a base.
We label the edge lengths as follows

\begin{center}
\begin{tabular}{cl}
\begin{minipage}{125pt}
\centerline{\includegraphics[width=125pt,keepaspectratio]{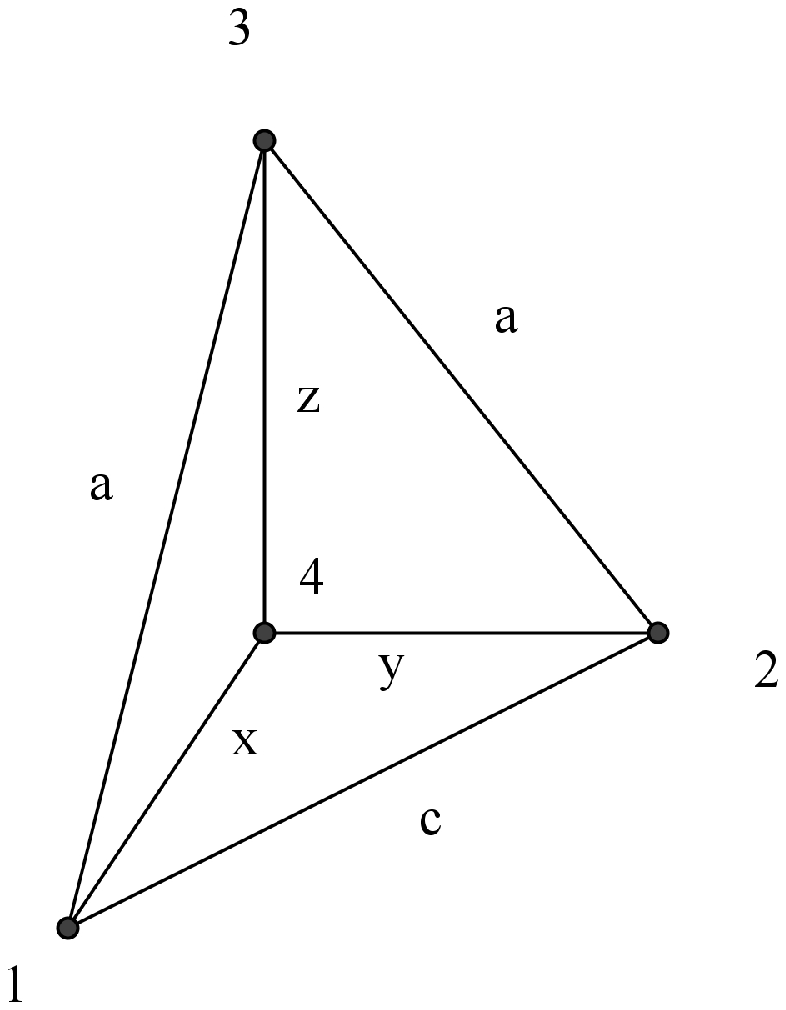}}
\end{minipage}
&
\begin{minipage}{50pt}
$r_{12}=c$\\
$r_{13}=b$\\
$r_{23}=a$\\
$r_{14}=x$\\
$r_{24}=y$\\
$r_{34}=z$
\end{minipage}
\end{tabular}
\end{center}

We have following obvious relations:
$a^2=y^2+z^2$,
$b^2=x^2+z^2$,
$c^2=x^2+y^2$.

By using them we can get
\begin{equation}
\label{E:trirec1}
\begin{array}{l@{\ =\ }l}
d_3(a,b,c) & 2(ax^2+by^2+cz^2-abc),\\
d_3(x,y,c) & 2xy(x+y-c),\\
d_3(x,b,z) & 2xz(x+z-b),\\
d_3(a,y,z) & 2yz(y+z-a)
\end{array}
\end{equation}
and
\begin{equation}
\label{E:trirec2}
\begin{array}{l}
{\Re(D_4)-64abcxyz-288V^2}=\\
\displaystyle
={4xyz}\sum_{cyc}2ax^2+\sum_{cyc}(2ab+cz+z^2)(x+y)-10abc
\end{array}
\end{equation}
where $\sum_{cyc}$ has three terms
\footnote{$\sum_{cyc}f(a,b,c,x,y,z)=f(a,b,c,x,y,z)+f(b,c,a,y,z,x)+f(c,a,b,z,x,y)$}
corresponding to a cycle $((a,x)\rightarrow (b,y)\rightarrow (c,z))$.

By writing $x+y=x+y-c+c$ and using the identity
\[
\sum_{cyc}c^2z=\sum_{cyc}(x^2+y^2)z=\sum_{cyc}(x+y)z^2=\sum_{cyc}z^2(x+y-c)+\sum_{cyc}ax^2
\]
we get that the second cyclic sum is equal to
\begin{equation}
\label{E:trirec3}
\sum_{cyc}(2ab+cz+z^2)(x+y)=6ab+
\sum_{cyc}(2ab+cz+2z^2)(x+y-c)+2\sum_{cyc}ax^2
\end{equation}
By inserting this into (\ref{E:trirec2}) we get
\[
\Re (D_4)-64abcxyz-288V^2=4xyz(2d_3(a,b,c)+\sum_{cyc} (2ab+cz+2z^2)(x+y-c))
\]
Hence $\Re (D_4)\geq 64abcxyz$ so the verification of the C2 of
Atiyah--Sutcliffe for trirectangular tetrahedra is finished.

%%%%%%%%%%%%%%%%%%%%%%%%%%%%%%%%%%%%%%%%%%%%%%%%%%%%%%%%%%%%%%%%%%%%%%%%%%%%%%
%
%     20060722
%
%%%%%%%%%%%%%%%%%%%%%%%%%%%%%%%%%%%%%%%%%%%%%%%%%%%%%%%%%%%%%%%%%%%%%%%%%%%%%%

\subsection{Atiyah--Sutcliffe conjectures for regular and semi--regular tetrahedra}

Semiregular (SR) tetrahedra are one of the simplest configurations of tetrahedra.
These tetrahedra have opposite edges equal and hence all faces are congruent.
Sometimes semi--regular tetrahedra are called isosceles tetrahedra.

\begin{center}
\begin{tabular}{cll}
\begin{minipage}{125pt}
\centerline{\includegraphics[width=125pt,keepaspectratio]{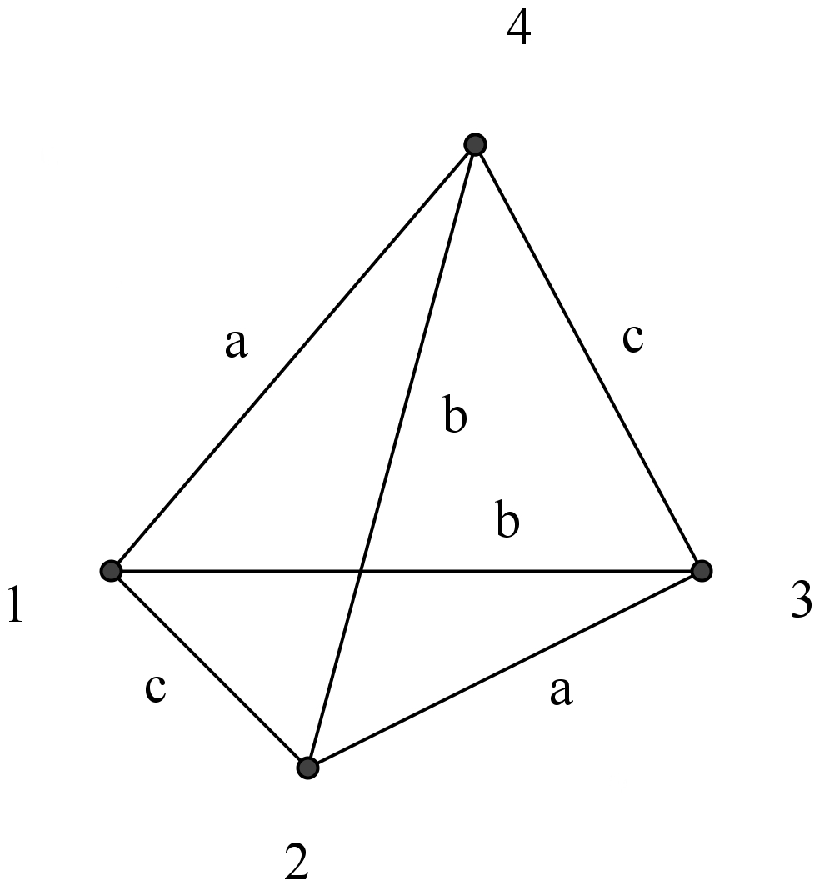}}
\end{minipage}
&
\begin{minipage}{150pt}
$r_{23}=r_{14}=a$\\
$r_{13}=r_{24}=b$\\
$r_{12}=r_{34}=c$
\end{minipage}
\end{tabular}
\end{center}

By (\ref{E:e0}) we get
\[
288V^2-4d_3(r_{12}r_{34},r_{13}r_{24},r_{14}r_{23})=0\
(\Rightarrow 288V^2=4d_3(r_{12}r_{34},r_{13}r_{24},r_{14}r_{23}))
\]
By (\ref{E:A4-a}) we get
\[
\begin{array}{l@{\ =\ }l}
A_4
& \displaystyle \sum_{l=1}^4(d_3(r_{il},r_{jl},r_{kl})+8r_{il}r_{jl}r_{kl})
                            d_3(r_{ij},r_{ik},r_{jk})\\[5mm]
& \displaystyle 4d_3(a,b,c)^2+32abc\, d_3(a,b,c)
\end{array}
\]
The quantity in the weak Four Points Conjecture is
\[
\begin{array}{l}
l.h.s-r.h.s=\\
\displaystyle\ =A_4-\left(4+\frac{3}{4}\right)288V^2-
 \sum_{l=1}^4\left(4+\frac{1}{4}\frac{d_3(r_{ij},r_{ik},r_{jk})}{r_{il}r_{jl}r_{kl}}\right)
 r_{ij}r_{ik}r_{jk}d_3(r_{ij},r_{ik},r_{jk})\\[5mm]
\displaystyle\ =4d_3(a,b,c)^2+32abc\, d_3(a,b,c)-(16+3)d_3(a^2,b^2,c^2)-[16abc\, d_3(a,b,c)+d_3(a,b,c)^2]\\[5mm]
\displaystyle\ =3(d_3(a,b,c)^2-d_3(a^2,b^2,c^2))+16(abc\, d_3(a,b,c)-d_3(a^2,b^2,c^2))\geq 0
\end{array}
\]
by using the inequalities $abc\geq d_3(a,b,c)$ and $d_3(a,b,c)^2\geq d_3(a^2,b^2,c^2)$
(see Appendix 2, Proposition \ref{P:A2}; also see \cite{Wolstenholme} or \cite{Walker}).

This proves the weak Four Points Conjecture for semiregular tetrahedra.\qed

The proof of the strong Four Points Conjecture for semiregular tetrahedra reduces to the positivity
of the following expression
\[
\displaystyle\ 4(d_3(a,b,c)^2-d_3(a^2,b^2,c^2))+15(abc\, d_3(a,b,c)-d_3(a^2,b^2,c^2))\geq 0
\]
which is also true by the same argument.

%%%%%%%%%%%%%%%%%%%%%%%%%%%%%%%%%%%%%%%%%%%%%%%%%%%%%%%%%%%%%%%%%%%%%%%%%%%%%%%
%
%    20060805
%
%%%%%%%%%%%%%%%%%%%%%%%%%%%%%%%%%%%%%%%%%%%%%%%%%%%%%%%%%%%%%%%%%%%%%%%%%%%%%%%

\subsection{Atiyah--Sutcliffe conjectures for parallelograms}

Given a parallelogram with vertices $1$, $2$, $3$ and $4$ denote by $a$, $b$ its side lengths and by
$e$, $f$ its diagonals.

\begin{center}
\begin{tabular}{cll}
\begin{minipage}{175pt}
\centerline{\includegraphics[width=130pt,keepaspectratio]{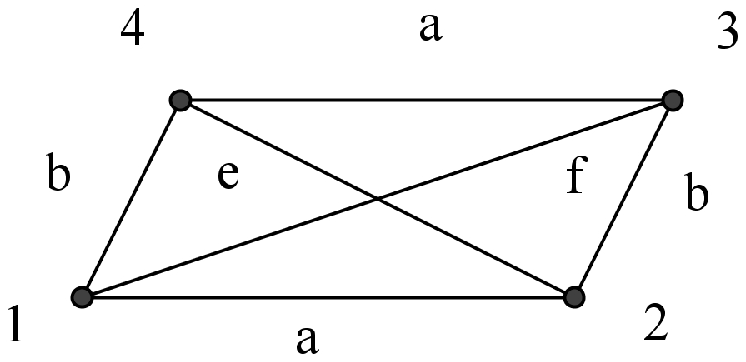}}
\end{minipage}
&
\begin{minipage}{150pt}
$r_{12}=r_{34}=a$\\
$r_{23}=r_{14}=b$\\
$r_{24}=e$\\
$r_{13}=f$
\end{minipage}
\end{tabular}
\end{center}

For the numbers $a,b,e,f$ we have the basic relation ("a parallelogram law")
\begin{equation}
\label{E:parallelogramlaw}
e^2+f^2=2(a^2+b^2)
\end{equation}
By using this relation we can rewrite various quantities in the Eastwood-Norbury formula.
\begin{proposition}
We have the following identities
\begin{enumerate}
  \item $d_3(a,b,e)=(a+b-e)(a-b+e)(-a+b+e)=(a+b-e)(a+b-f)(a+b+f)$
  \item $\Delta:=(a+b+e)d_3(a,b,e)=(a+b+f)d_3(a,b,f)=$\\
        $=(a+b+e)(a+b+f)(a+b-e)(a+b-f)=$\\
        $=2a^2b^2+2a^2e^2+2b^2e^2-a^4-b^4-e^4=2a^2b^2+2a^2f^2+2b^2f^2-a^4-b^4-f^4$
  \item $4ab+e^2-f^2=2(a+b+f)(a+b-f)$, $4ab+f^2-e^2=2(a+b+e)(a+b-e)$
  \item $d_3(a^2,b^2,ef)=(a^2+b^2-ef)\Delta$
  \item $d_3(a,b,e)d_3(a,b,f)-d_3(a^2,b^2,ef)=(2ab-2ef-(a+b)(e+f))\Delta$
  \item $ed_3(a,b,f)+fd_3(a,b,e)=(a+b-e)(a+b-f)(e^2+f^2+(a+b)(e+f))$
  \item $(4ab+e^2-f^2)ed_3(a,b,f)+(4ab+f^2-e^2)fd_3(a,b,e)=2((a+b)(e+f)-2ef)\Delta$
\end{enumerate}
\end{proposition}
\proof
For {\it 1.} we write $(a-b+e)(-a+b+e)=e^2-(a-b)^2=2a^2+2b^2-f^2-(a-b)^2=(a+b-f)(a+b+f)$.
Identity {\it 2.} follows from {\it 1.} directly. For {\it 3.} we substitute
$e^2=2a^2+2b^2-f^2$ and simplify. For {\it 4.} we compute and use {\it 2.}:
\[
(a^2-b^2+ef)(-a^2+b^2+ef)=e^2f^2-(a^2-b^2)^2=e^2(2a^2+2b^2-e^2)+2a^2b^2-a^4-b^4=\Delta
\]
For {\it 5.} we first use {\it 1.} and then {\it 4.}:
$d_3(a,b,e)d_3(a,b,f)-d_3(a^2,b^2,ef)$ $=$ $(a+b+f)(a+b-e)(a+b-f)d_3(a,b,f)$ $-$
$(a^2+b^2-ef)\Delta$ $=$ $[(a+b)^2-(a+b)(e+f)+ef]\Delta$ $-$ $(a^2+b^2-ef)\Delta$
$=$ $[2ab+2ef-(a+b)(e+f)]\Delta$

For {\it 6.} we use {\it 1.} twice.

For {\it 7.} we first use {\it 3.} and then {\it 2.}:
\[
\begin{array}{l}
l.h.s.=2(a+b+f)(a+b-f)ed_3(a,b,f)+2(a+b+e)(a+b-e)fd_(a,b,e)\\
\hphantom{l.h.s.}=2[(a+b-f)e+(a+b-e)f]\Delta
\end{array}
\]\qed

Now we apply Eastwood-Norbury formula (note that $288V^2=0$, $D_4=real$)
\[
D_4-64\prod r_{ij}=-4d_3(a^2,b^2,c^2)+A_4
\]
where
\[
\begin{array}{l@{\ =\ }l}
A_4
&
2[d_3(a,b,e)+8abc+e(e^2-f^2)]d_3(a,b,f)+
2[d_3(a,b,f)+8abc+f(f^2-e^2)]d_3(a,b,e)\\
&
I_0+I_1+I_2
\end{array}
\]
where
\[
\begin{array}{l@{\ =\ }l}
I_0 & 4d_3(a,b,e)d_3(a,b,f)\\
I_1 & 2[4abe+e(e^2-f^2)]d_3(a,b,f)+2[4abf+f(f^2-e^2)]d_3(a,b,e)\\
    & 4((a+b)(e+f)-2ef)\Delta\ \ (\mbox{by {7.}})\\
I_2 & 2[4abe\,d_3(a,b,f)+4abf\,d_3(a,b,e)]\\
    & 8ab(a+b-e)(a+b-f)(e^2+f^2+(a+b)(e+f))\ \ (\mbox{by {6.}})
\end{array}
\]
By using {\it 5.} we have
\[
\begin{array}{l@{\ =\ }l}
D_4-64\prod r_{ij}
& 4(d_3(a,b,e)d_3(a,b,f)-d_3(a^2,b^2,ef))+I_1+I_2\\
& 4((2ab+2ef-(a+b)(e+f))\Delta+((a+b)(e+f)-2ef)\Delta)+I_2\\
& 8ab\Delta+I_2\geq 0
\end{array}
\]
This proves the Atiyah--Sutcliffe conjecture (C2) for parallelograms.
The Atiyah--Sutcliffe conjecture (C3) for parallelograms
\[
D_4^2\geq (d_3(a,b,e)+8abe)^2(d_3(a,b,f)+8abf)^2
\]
is equivalent to  the positivity of
\[
D_4-d_3(a,b,e)d_3(a,b,f)-8[abf\, d_3(a,b,e)+abe\, d_3(a,b,f)]-64a^2b^2ef\geq 0\\
\]
but we can prove even stronger statement
\[
\begin{array}{l}
D_4-2d_3(a,b,e)d_3(a,b,f)-8[abf\, d_3(a,b,e)+abe\, d_3(a,b,f)]-64a^2b^2ef=\\
=8ab\Delta-2d_3(a,b,e)d_3(a,b,f)%\\
=2[4ab-(a+b-e)(a+b-f)]\Delta\geq 0
\end{array}
\]
because the triangle inequalities $b<e+a$ and $a<f+b$ imply
\[
(a+b-e)(a+b-f)<2a\cdot 2b=4ab.
\]
Thus we have verified also C3 for parallelograms.

Finally we verify our strong Four Point Conjecture for parallelograms as follows
\[
\begin{array}{l}
D_4-64\prod r_{ij}-\sum(4+\frac{1}{4})r_{il}r_{jl}r_{kl}d_3(r_{ij},r_{jl},r_{ik})\\
=8ab\Delta -\frac{1}{4}(I_2/4)\\
=8ab(a+b-e)(a+b-f)[(a+b+e)(a+b+f)-\frac{1}{16}(e^2+f^2+(a+b)(e+f))]\\
=\frac{1}{2}ab(a+b-e)(a+b-f)[16((a+b)^2+(a+b)(e+f)+ef)-(2a^2+2b^2+(a+b)(e+f))]\\
=\frac{1}{2}ab(a+b-e)(a+b-f)[14(a^2+b^2)+32ab+15(a+b)(e+f)+16ef]\geq 0
\end{array}
\]\qed

%%%%%%%%%%%%%%%%%%%%%%%%%%%%%%%%%%%%%%%%%%%%%%%%%%%%%%%%%%%%%%%%%%%%%%%%%%%%%%%
%
%  20060808
%
%%%%%%%%%%%%%%%%%%%%%%%%%%%%%%%%%%%%%%%%%%%%%%%%%%%%%%%%%%%%%%%%%%%%%%%%%%%%%%%

\subsection{Atiyah--Sutcliffe conjectures for "wedge" tetrahedra}

A tetrahedron with two pairs of opposite edges having the same length we simply call
a "wedge" tetrahedron.

\begin{center}
\begin{tabular}{cll}
\begin{minipage}{125pt}
\centerline{\includegraphics[width=100pt,keepaspectratio]{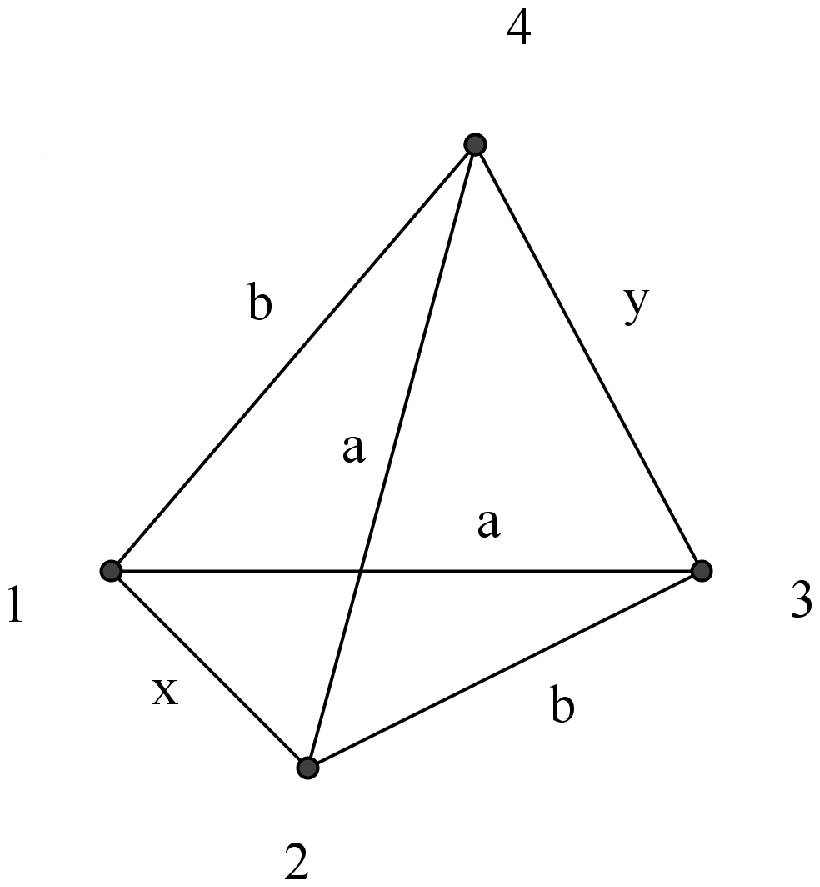}}
\end{minipage}
&
\begin{minipage}{150pt}
$r_{12}=x$\\
$r_{34}=y$\\
$r_{13}=r_{24}=a$\\
$r_{23}=r_{14}=b$
\end{minipage}
\end{tabular}
\end{center}

\begin{center}
\begin{tabular}{cll}
\begin{minipage}{150pt}
If $x=y=c$ we get a semiregular tetrahedron and if all points lie in a plane then we get either
a parallelogram or an isosceles trapezium.
\end{minipage}
&
\begin{minipage}{200pt}
\centerline{\includegraphics[width=200pt,keepaspectratio]{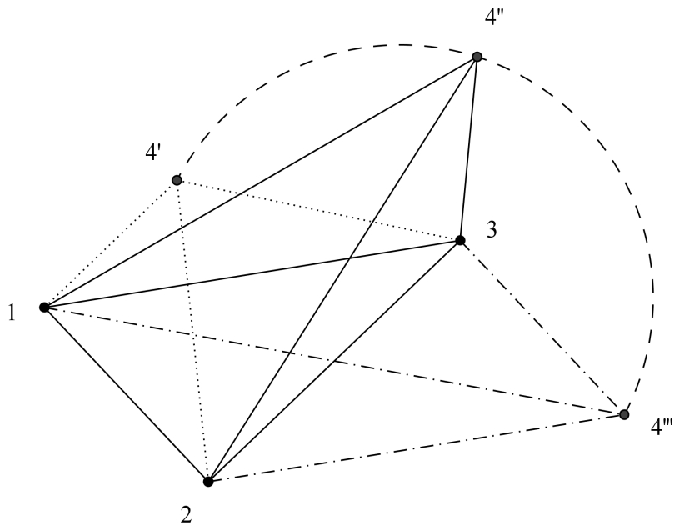}}
\end{minipage}
\end{tabular}
\end{center}

Again we compute the data appearing in the Eastwood--Norbury formula
\begin{equation}
\label{E:wedge2}
\begin{array}{l}
-4d_3(r_{12}r_{34},r_{13}r_{24},r_{14}r_{23})=\\
=-4d_3(xy,a^2,b^2)\\
=-4(xy-a^2+b^2)(xy+a^2-b^2)(a^2+b^2-xy)
\end{array}
\end{equation}
and we have the basic inequalities
\begin{equation}
\label{E:wedge3}
xy+b^2\geq a^2,\ xy+a^2\geq b^2,\ a^2+b^2\geq xy
\end{equation}
The positivity of the volume
\begin{equation}
\label{E:wedge4}
144V^2=xy-a^2+b^2)(xy+a^2-b^2)(2a^2+2b^2-x^2-y^2)
\end{equation}
gives us one more basic inequality
\begin{equation}
\label{E:wedge5}
2a^2+2b^2\geq x^2+y^2
\end{equation}
We have
\begin{equation}
\label{E:wedge6}
\begin{array}{ll}
A_4=\!\!\!\!
&2[a((b+x)^2-a^2)+b((a+x)^2-b^2)+x((a+b)^2-y^2)]d_3(a,b,y)\\
&+2[a((b+y)^2-a^2)+b((a+y)^2-b^2)+y((a+b)^2-x^2)]d_3(a,b,x)
\end{array}
\end{equation}
By using identity
\begin{equation}
\label{E:wedge7}
d_3(a,b,c)=a(b^2+c^2-a^2)+b(a^2+c^2-b^2)+c(a^2+b^2-c^2)-2abc
\end{equation}
we can rewrite $A_4$ as follows
\begin{equation}
\label{E:wedge8}
\begin{array}{ll}
A_4=\!\!\!\!
&2[4abx+d_3(a,b,x)-x(a^2+b^2-x^2)+2abx+x((a+b)^2-y^2)]d_3(a,b,y)\\
&+2[4aby+d_3(a,b,y)-y(a^2+b^2-y^2)+2aby+y((a+b)^2-x^2)]d_3(a,b,x)\\
\hphantom{A_4}=\!\!\!\!
&2[4abx+d_3(a,b,x)-x((a-b)^2-x^2)+x((a+b)^2-y^2)]d_3(a,b,y)\\
&+2[4aby+d_3(a,b,y)-y((a-b)^2-y^2)+y((a+b)^2-x^2)]d_3(a,b,x)\\
\hphantom{A_4}=\!\!\!\!
&\{8abx+d_3(a,b,x)+[d_3(a,b,x)-2x((a-b)^2-x^2)]+2x((a+b)^2-y^2)\}d_3(a,b,y)\\
&+\{8aby+d_3(a,b,y)+[d_3(a,b,y)-2y((a-b)^2-y^2)]+2y((a+b)^2-x^2)\}d_3(a,b,x)\\
\end{array}
\end{equation}
Now we compute
\[
\begin{array}{l}
d_3(a,b,x)-2x((a-b)^2-x^2)=\\
=(a+b-x)(a-b+x)(-a+b+x)+2x(a-b+x)(-a+b+x)\\
=(a+b+x)(a-b+x)(-a+b+x)
\end{array}
\]
The contribution $A_4^{[\ ]}$ of both square brackets in $A_4$ is equal to
\begin{equation}
\begin{array}{l}
A_4^{[\ ]}:=[d_3(a,b,x)-2x((a-b)^2-x^2)]d_3(a,b,y)+\\
\hphantom{A_4^{[\ ]}:=}+[d_3(a,b,y)-2y((a-b)^2-y^2)]d_3(a,b,x)\\
\hphantom{A_4^{[\ ]}:}\! =
(x^2-(a-b)^2)(y^2-(a-b)^2)[(a+b+x)(a+b-y)+(a+b+y)(a+b-x)]\\
\hphantom{A_4^{[\ ]}:}\! =
(x^2-(a-b)^2)(y^2-(a-b)^2)(2(a+b)^2-2xy)\\
\hphantom{A_4^{[\ ]}:}\! =
4ab(x^2-(a-b)^2)(y^2-(a-b)^2)+2(x^2-(a-b)^2)(y^2-(a-b)^2)(a^2+b^2-xy)
\end{array}
\end{equation}
At this point we have discovered the following beautiful identity
\begin{equation}
\label{E:wedge*}
\begin{array}{l}
[x^2-(a-b)^2][y^2-(a-b)^2]=\\
=(xy-a^2+b^2)(xy+a^2-b^2)+(a-b)^2(2a^2+2b^2-x^2-y^2)
\end{array}
\end{equation}
By this identity we can write
\[
\begin{array}{l}
A_4^{[\ ]}=4ab(x^2-(a-b)^2)(y^2-(a-b)^2)\\
\hphantom{A_4^{[\ ]}=}+2(a-b)^2(2a^2+2b^2-x^2-y^2)(a^2+b^2-xy)\\
\hphantom{A_4^{[\ ]}=}+2d_3(a^2,b^2,xy)
\end{array}
\]
\begin{lemma}
\label{L:2.10}
We have the following inequality for "wedge" tetrahedra
\[
d_3(a^2,b^2,xy)\leq 2ab(x^2-(a-b)^2)(y^2-(a-b)^2)
\]
\end{lemma}
\proof
Recall that
\[
d_3(a^2,b^2,xy)=(a^2+b^2-xy)(a^2-b^2+xy)(-a^2+b^2+xy)
\]
Let $a\geq b$. Then the triangle inequalities $a\leq b+x$ and $a\leq b+y$ imply
$(a-b)^2\leq xy$ i.e.\ $a^2+b^2-xy\leq 2ab$.
Since $2a^2+2b^2-x^2-y^2\geq 0$ (inequality (\ref{E:wedge5})) then from our inequality
(\ref{E:wedge*}) it follows that
\[
(a^2-2b^2+xy)(-a^2+b^2+xy)\leq (x^2-(a-b)^2)(y^2-(a-b)^2)
\]
By multiplying the last two inequalities Lemma follows.\qed

As a consequence of Lemma we get immediately that
\[
A_4\geq A_4^{[\ ]}\geq 4d_3(a^2,b^2,xy)
\]
because the remaining terms in $A_4$ are all nonnegative.
This verifies the A--S conjecture C2 for "wedge" tetrahedra.
\begin{remark}
Instead of splitting $2(a+b)^2-2xy=4ab+2(a^2+b^2-xy)$ (used above),
we can use the identity
\[
2(a+b)^2-2xy=4(a^2+b^2-xy)+2(xy-(a-b)^2)
\]
to obtain explicit formula for $A_4^{[\ ]}$:
\[
\begin{array}{l}
A_4^{[\ ]}=\\
{[}4(a^2+b^2-xy)+2(xy-(a-b)^2)][(xy-a^2+b^2)(xy+a^2-b^2)+\\
+(a-b)^2(2a^2+2b^2-x^2-y^2)]=\\
=4d_3(a^2,b^2,xy)+4(a^2+b^2-xy)(2a^2+2b^2-x^2-y^2)(a-b)^2+\\
+2(xy-(a-b)^2)(x^2-(a-b)^2)(y^2-(a-b)^2)
\end{array}
\]
which, without using Lemma \ref{L:2.10}, implies inequality
\[
A_4^{[\ ]}\geq 4d_3(a^2,b^2,xy)
\]
needed for the verification of A--S conjecture C2 for "wedge" tetrahedra.
\end{remark}

Now we state a final formula for "wedge" tetrahedra:
\[
\begin{array}{|ll|}\hline
&\mbox{\underline{First explicit formula for wedge tetrahedra}:}\\[2mm]
\Re(D_4)=
&(d_3(a,b,x)+8abx)(d_3(a,b,y)+8aby)+d_3(a,b,x)d_3(a,b,y)+\\
&+2x((a+b)^2-y^2)d_3(a,b,y)+2y((a+b)^2-x^2)d_3(a,b,x)+\\
&+4(a^2+b^2-xy)(2a^2+2b^2-x^2-y^2)(a-b)^2+\\
&+2(xy-(a-b)^2)(x^2-(a-b)^2)(y^2-(a-b)^2)\\
&+288V^2\\\hline
\end{array}
\]
which implies a strengthened A--S conjecture C3 for wedge tetrahedra
\[
\begin{array}{rl}
\Re(D_4)\geq
&(d_3(a,b,x)+8abx)(d_3(a,b,y)+8aby)+d_3(a,b,x)d_3(a,b,y)+288V^2\\
\geq
&(d_3(a,b,x)+8abx)(d_3(a,b,y)+8aby)
\end{array}
\]
In the sequel we obtain an alternative formula for the real part of the Atiyah
determinant for a wedge tetrahedra.

We group terms in $A_4$ differently as follows:
\[
\begin{array}{rl}
A_4=
&2[4abx+d_3(a,b,x)+x(4ab+x^2-y^2)]d_3(a,b,y)+\\
&+2[4aby+d_3(a,b,y)+x(4ab+y^2-x^2)]d_3(a,b,x)
\end{array}
\]
By letting
\[
2s^2+2b^2-x^2-y^2=:2h\ \ (\geq 0)
\]
we can rewrite
\[
4ab+x^2-y^2=4ab+x^2+(2h+x^2-2a^2-2b^2)=2(h+x^2-(a-b)^2)
\]
and similarly for
\[
4ab+y^2-x^2=2(h+y^2-(a-b)^2)
\]
Thus
\[
\begin{array}{rl}
A_4=
&4d_3(a,b,x)d_3(a,b,y)+8abx\, d_3(a,b,y)+8aby\, d_3(a,b,x)+\\
&+4h(x\, d_3(a,b,y)+y\, d_3(a,b,x))+4A'_4
\end{array}
\]
where
\begin{equation}
\label{E:wedgeo}
\begin{array}{rl}
A_4'
&=x(x^2-(a-b)^2)d_3(a,b,y)+y(y^2-(a-b)^2)d_3(a,b,x)\\
&=(x^2-(a-b)^2)(y^2-(a-b)^2)[x(a+b-y)+y(a+b-x)]\\
&=(x^2-(a-b)^2)(y^2-(a-b)^2)[(x-y)^2+x(a+b-x)+y(a+b-y)]\\
&=[(xy-a^2+b^2)(xy+a^2-b^2)+2(a-b)^2h][(x-y)^2+x(a+b-x)+y(a+b-y)]\\
\end{array}
\end{equation}
by our identity (\ref{E:wedge*}).

Note that
\[
-144V^2+2d_3(a^2,b^2,xy)=(xy-a^2+b^2)(xy+a^2-b^2)(x-y)^2
\]
So
\[
\begin{array}{l}
A_4'=(2d_3(a^2,b^2,xy)-144V^2)+2(a-b)^2h[x(a+b-y)+y(a+b-x)]\\
\hphantom{A_4'=}+(xy-a^2+b^2)(xy+a^2-b^2)(x(a+b-x)+y(a+b-y))
\end{array}
\]
By writing
\[
\begin{array}{l}
4A_4'=2A_4'+2A_4'=\\
=2(x^2-(a-b)^2)(y^2-(a-b)^2)[x(a+b-y)+y(a+b-x)]+\\
+\{4d_3(a^2,b^2,xy)-288V^2+4(a-b)^2h[x(a+b-y)+y(a+b-x)]+\\
+2(xy-a^2+b^2)(xy+a^2-b^2)(x(a+b-x)+y(a+b-y))\}=\\
=4d_3(a^2,b^2,xy)-288V^2+[2(x^2-(a-b)^2)(y^2-(a-b)^2)+4(a-b)^2h]\cdot\\
\cdot (x(a+b-y)+y(a+b-x))+\\+
2(xy-a^2+b^2)(xy+a^2-b^2)[x(a+b-x)+y(a+b-y)]
\end{array}
\]
we obtain the following explicit formula for the real part of Atiyah determinant for
"wedge" tetrahedron:
\[
\begin{array}{|ll|}\hline
&\mbox{\underline{Second explicit formula for wedge tetrahedra}:}\\[2mm]
\Re(D_4)=\!\!\!\!
&(d_3(a,b,x)+8abx)(d_3(a,b,y)+8aby)+3d_3(a,b,x)d_3(a,b,y)+\\
&+2x((a+b)^2-y^2)d_3(a,b,y)+2y((a+b)^2-x^2)d_3(a,b,x)+\\
&2(x^2y^2-(a^2-b^2))[x(a+b-x)+y(a+b-y)]+\\
&+2[(a-b)^2(x(a+b-y)+y(a+b-x))](2a^2+2b^2-x^2-y^2)\\\hline
\end{array}
\]
which implies another strengthening of the Atiyah--Sutcliffe conjecture C3 for "wedge"
tetrahedra
\[
\begin{array}{l}
\Re (D_4)\geq (d_3(a,b,x)+8abx)(d_3(a,b,y)+8aby)+3d_3(a,b,x)d_3(a,b,y)\\
\hphantom{\Re (D_4)}\geq (d_3(a,b,x)+8abx)(d_3(a,b,y)+8aby)
\end{array}
\]

%%%%%%%%%%%%%%%%%%%%%%%%%%%%%%%%%%%%%%%%%%%%%%%%%%%%%%%%%%%%%%%%%%%%%%%%%%%%%%%

\subsection{Atiyah determinant for triangles and quadrilaterals via trigonometry}
\label{S:Adftaqvt}
Denote the three points $x_1$, $x_2$, $x_3$ simply by symbols $1, 2, 3$ and let $X$, $Y$ and $Z$ denote
the angles of the triangle at vertices $1$, $2$ and $3$ respectively. Then we can express
the Atiyah determinant $D_3=d_3(r_{12},r_{13},r_{23})+8r_{12}r_{13}r_{23}$ as follows
\[
D_3=4r_{12}r_{13}r_{23}\left(\cos^2\frac{X}{2}+\cos^2\frac{Y}{2}+\cos^2\frac{Z}{2}\right).
\]
This follows, by using cosine law and sum to product formula for cosine, from the following identity
\[
\begin{array}{l@{\ =\ }l}
d_3(a,b,c)+8abc&(a+b-c)(a-b+c)(-a+b+c)+8abc\\
               &a((b+c)^2-a^2)+b((c+a)^2-b^2)+c((a+b)^2-c^2).
\end{array}
\]
Now we shall translate the Eastwood--Norbury formula for (planar quadrilaterals) into a trigonometric form.
Denote the four points $x_1$, $x_2$, $x_3$, $x_4$ simply by symbols $1, 2, 3, 4$ and denote by
\[
(X^{(1)}, Y^{(1)}, Z^{(1)}),\ \
(X^{(2)}, Y^{(2)}, Z^{(2)}),\ \
(X^{(3)}, Y^{(3)}, Z^{(3)}),\ \
(X^{(4)}, Y^{(4)}, Z^{(4)})
\]
 the angles of the triangles $234$, $341$, $412$, $123$ in this cyclic order
 (i.e.\ the angle of a triangle $412$ at vertex $2$ is $Z^{(3)}$ etc.).

Next we denote by $c_l$, ($1\leq l\leq 4$), the sums of cosines squared of half-angles of the $l$--th triangle i.e.:
\[
c_l:=\cos^2\frac{X^{(l)}}{2}+\cos^2\frac{Y^{(l)}}{2}+\cos^2\frac{Z^{(l)}}{2},\ \ l=1,2,3,4.
\]
Similarly, we denote by $\widehat{c}_l$, ($1\leq l\leq 4$),
the sum of cosines squared of half-angles at the $l$--th vertex of our quadrilateral thus
\[
\begin{array}{c}
\displaystyle \widehat{c}_1=\cos^2\frac{Z^{(2)}}{2}+\cos^2\frac{Y^{(3)}}{2}+\cos^2\frac{X^{(4)}}{2}\\[3mm]
\displaystyle \widehat{c}_2=\cos^2\frac{Z^{(3)}}{2}+\cos^2\frac{Y^{(4)}}{2}+\cos^2\frac{X^{(1)}}{2}\\[3mm]
\displaystyle \widehat{c}_3=\cos^2\frac{Z^{(4)}}{2}+\cos^2\frac{Y^{(1)}}{2}+\cos^2\frac{X^{(2)}}{2}\\[3mm]
\displaystyle \widehat{c}_4=\cos^2\frac{Z^{(1)}}{2}+\cos^2\frac{Y^{(2)}}{2}+\cos^2\frac{X^{(3)}}{2}
\end{array}
\]
 Then the term
$A_4$ in the Eastwood--Norbury formula can be rewritten as
\[
\begin{array}{r@{\ =\ }l}
  A_4 & \displaystyle\sum_{l=1}^4(4r_{li}r_{lj}r_{lk}\widehat{c}_l)\cdot 4r_{ij}r_{ik}r_{jk}(c_l-2) \\
    & \displaystyle16r_{12}r_{13}r_{23}r_{14}r_{24}r_{34}\sum_{l=1}^4\widehat{c}_l(c_l-2). \\
\end{array}
\]
where for each $l$ we write $\{1,2,3,4\}\setminus\{l\}=\{i<j<k\}$.\\

In order to rewrite the term $-4d_3(r_{12}r_{34},r_{13}r_{24},r_{14}r_{23})$ into a trigonometric form
we recall a theorem
of M\" obius (\cite{Moebius}) which claims that for any quadrilateral $1234$ in a plane the products
$r_{12}r_{34}$, $r_{13}r_{24}$ and $r_{14}r_{23}$ are proportional to the sides of a triangle
whose angles are the differences of angles in the quadrilateral $1234$:
\[
\begin{array}{c@{\ =\ }c}
  X & \sphericalangle 134 - \sphericalangle 124\\
  Y & \sphericalangle 214 - \sphericalangle 234\\
  Z & \sphericalangle 413 - \sphericalangle 423
%  X & \sphericalangle 243 - \sphericalangle 213 & X^{(4)}-Z^{(1)} \\
%  Y & \sphericalangle 341 - \sphericalangle 321 & Y^{(2)}-(-Y^{(4)}) \\
%  Z & \sphericalangle 142 - \sphericalangle 132 & -X^{(3)}+Z^{(4)}
\end{array}
\]

%%%%%%%%%%%%%%%%%%%%%%%%%%%%%%%%%%%%%%%%%%%%%%%%%%%%%%%%%%%%%%%%%%%%%%%%%%%%%%% [8]
Thus
\[
-4d_3(r_{12}r_{34},r_{13}r_{24},r_{14}r_{23})=-16r_{12}r_{13}r_{23}r_{14}r_{24}r_{34}(c-2)
\]
where
\[
c=\cos^2\frac{X}{2}+\cos^2\frac{Y}{2}+\cos^2\frac{Z}{2}.
\]
Thus we have obtained a trigonometric formula for Atiyah determinant of quadrilaterals
\[
\begin{array}{l@{\ =\ }l}
\Re(D_4)&\displaystyle\prod_{1\leq i<j\leq 4}r_{ij}\left(64-16(c-2)+16\sum_{l=1}^4\widehat{c}_l(c_l-2)\right)\\[4mm]
&\displaystyle 16\prod_{1\leq i<j\leq 4}r_{ij}\left(6-c+\sum_{l=1}^4\widehat{c}_l(c_l-2)\right)
\end{array}
\]
Now we shall verify Atiyah--Sutcliffe conjecture for cyclic quadrilaterals.

\begin{center}
\begin{tabular}{cl}
\begin{minipage}{100pt}
\centerline{\includegraphics[width=100pt,keepaspectratio]{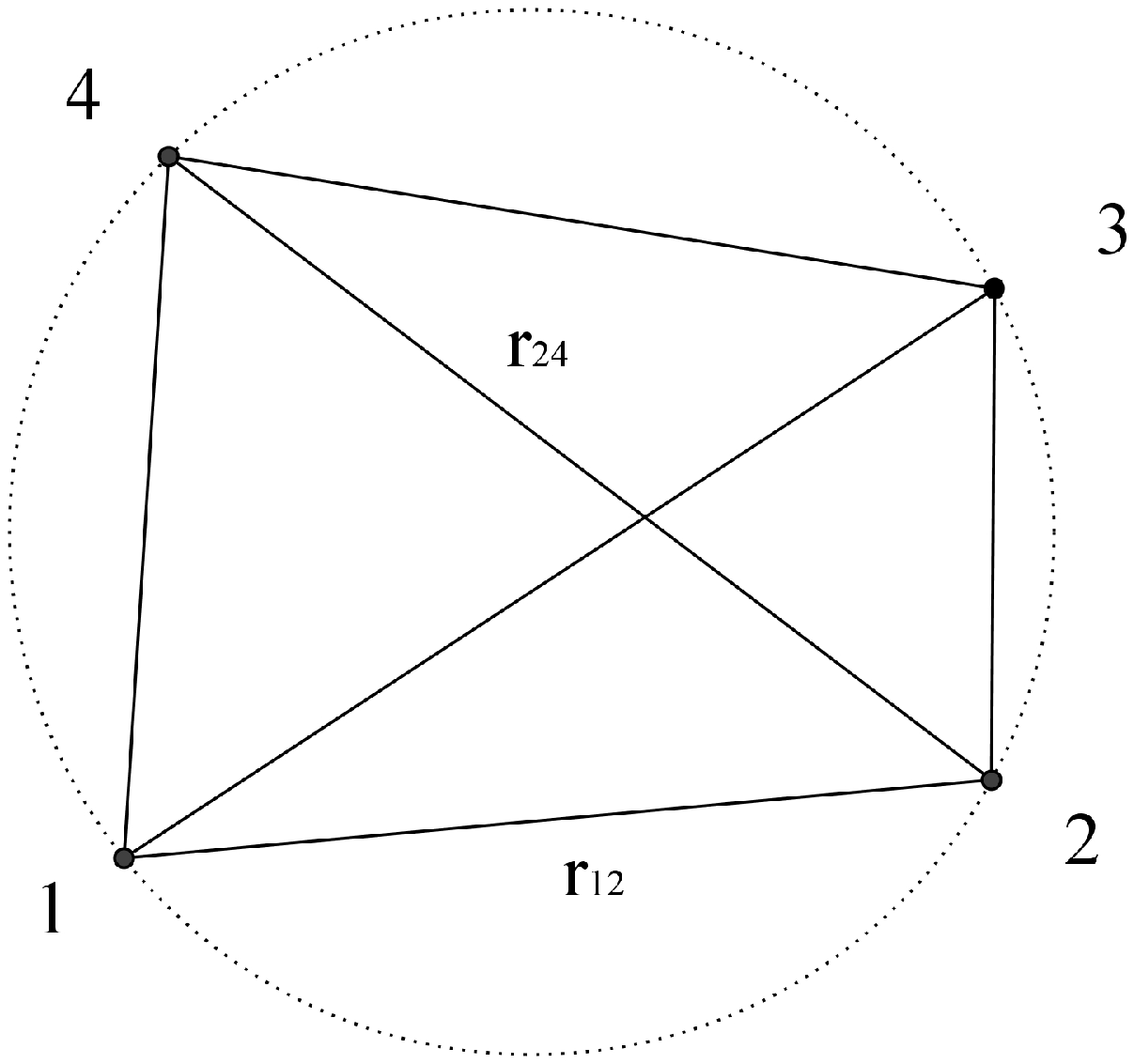}}
\end{minipage}
&
\begin{minipage}{120pt}
\begin{center}
Ptolemy's theorem\\[3mm]
$r_{12}r_{34}+r_{23}r_{14}=r_{13}r_{24}$
\end{center}
\end{minipage}
\end{tabular}
\end{center}
In this case, by a well known Ptolemy's theorem, we see that
\[
-4d_3(r_{12}r_{34},r_{13}r_{24},r_{14}r_{23})=0\ \ (\Leftrightarrow c=2)
\]
By using the equality of angles $Z^{(2)}=X^{(1)}$, $Z^{(3)}=X^{(2)}$, $Z^{(4)}=X^{(3)}$, $Z^{(1)}=X^{(4)}$
and $Y^{(1)}+Y^{(3)}=\pi=Y^{(2)}+Y^{(4)}$ (angles with vertex on a circle's circumference with the same
endpoints are equal or supplement of each other)we obtain
\[
\begin{array}{c}
\displaystyle \widehat{c}_1=\cos^2\frac{X^{(1)}}{2}+\sin^2\frac{Y^{(1)}}{2}+\cos^2\frac{Z^{(1)}}{2}=c_1-\cos Y^{(1)},\\[3mm]
\displaystyle \widehat{c}_2=\cos^2\frac{X^{(2)}}{2}+\sin^2\frac{Y^{(2)}}{2}+\cos^2\frac{Z^{(2)}}{2}=c_2-\cos Y^{(2)},\\[3mm]
\displaystyle \widehat{c}_3=\cos^2\frac{X^{(3)}}{2}+\sin^2\frac{Y^{(3)}}{2}+\cos^2\frac{Z^{(3)}}{2}=c_3-\cos Y^{(3)},\\[3mm]
\displaystyle \widehat{c}_4=\cos^2\frac{X^{(4)}}{2}+\sin^2\frac{Y^{(4)}}{2}+\cos^2\frac{Z^{(4)}}{2}=c_4-\cos Y^{(4)}.
\end{array}
\]
Now we have
\[
\begin{array}{rl}
\Re (D_4)=&
    \displaystyle \left(\prod_{1\leq i<j\leq 4}r_{ij} \right)
    \left(64+16\sum_{l=1}^4(c_l-\cos Y^{(l)})(c_l-2)\right)\\[6mm]
\geq&
\displaystyle \left(\prod_{1\leq i<j\leq 4}r_{ij} \right)
    \left(64+16\sum_{l=1}^4(c_l-1)(c_l-2)\right)\\[6mm]
\end{array}
\]
(here we have used that $2\leq c_l (\leq\frac{9}{4})$ for each $l=1,2,3,4$)
\[
\begin{array}{rl}
\hphantom{\Re (D_4)}\geq&
\displaystyle \left(\prod_{1\leq i<j\leq 4}r_{ij} \right)
    \left(64+16\sum_{l=1}^4(c_l-2)+16\sum_{l=1}^4(c_l-2)^2\right)\\[6mm]
\geq&
\displaystyle \left(\prod_{1\leq i<j\leq 4}r_{ij} \right)
    \left(64+16\sum_{l=1}^4(c_l-2)+4\left(\sum_{l=1}^4(c_l-2)\right)^2\right)
\end{array}
\]
(by quadratic--arithmetic inequality)
\[
\begin{array}{rl}
\hphantom{\Re (D_4)}=&
\displaystyle \left(\prod_{1\leq i<j\leq 4}r_{ij} \right)
    \left(\left(8+\sum_{l=1}^4(c_l-2)\right)^2+3\left(\sum_{l=1}^4(c_l-2)\right)^2\right)\\[6mm]
=&
\displaystyle \left(\prod_{1\leq i<j\leq 4}r_{ij} \right)
    \left(\left(\sum_{l=1}^4c_l\right)^2+3\left(\sum_{l=1}^4(c_l-2)\right)^2\right)\\[6mm]
\geq&
\displaystyle \left(\prod_{1\leq i<j\leq 4}r_{ij} \right)
    \left(\sum_{l=1}^4c_l\right)^2\geq 16\sqrt{c_1c_2c_3c_4}\prod_{1\leq i<j\leq 4}r_{ij}
\end{array}
\]
by A--G inequality.

Finally,
\[
\begin{array}{rl}
|D_4|^2=&\displaystyle
|\Re(D_4)|^2\geq 4^4c_1c_2c_3c_4\prod_{1\leq i<j\leq 4}r_{ij}^2\\
=&\displaystyle
    \prod_{l=1}^4(4r_{ij}r_{ik}r_{jk}c_l)
=
    \prod_{l=1}^4(d_3(r_{ij},r_{ik},r_{jk})+8r_{ij}r_{ik}r_{jk})
\end{array}
\]
where for each $l$ we write $\{1,2,3,4\}\setminus\{l\}=\{i<j<k\}$.
This finishes verification of Atiyah--Sutcliffe conjectures for cyclic quadrilaterals.

\begin{proposition}
The weak Four Points Conjecture for cyclic quadrilaterals holds true.
\end{proposition}
\proof
From the formula obtained above we proceed along a different path

\[
\begin{array}{rl}
\Re (D_4)=&
    \displaystyle \left(\prod_{1\leq i<j\leq 4}r_{ij} \right)
    \left(64+16\sum_{l=1}^4(c_l-\cos Y^{(l)})(c_l-2)\right)\\[6mm]
\hphantom{\Re (D_4)}=&
\displaystyle \left(\prod_{1\leq i<j\leq 4}r_{ij} \right)
    \left(4\sum_{l=1}^4\left[4+4(c_l-\cos Y^{(l)})(c_l-2)\right]\right)\\[6mm]
%\hphantom{\Re (D_4)}=&
%\displaystyle \left(\prod_{1\leq i<j\leq 4}r_{ij} \right)
%    \left(4\sum_{l=1}^4\left[c_l^2+(c_l-2)[4(c_l-\cos Y^{(l)})-2-c_l]\right]\right)\\[6mm]
\hphantom{\Re (D_4)}=&
\displaystyle \left(\prod_{1\leq i<j\leq 4}r_{ij} \right)
    \left(4\sum_{l=1}^4\left[c_l^2+(c_l-2)[3(c_l-2) +4(1-\cos Y^{(l)})]\right]\right)\\[6mm]
\hphantom{\Re (D_4)}\geq&
\displaystyle \left(\prod_{1\leq i<j\leq 4}r_{ij} \right)
    \left(4\sum_{l=1}^4c_l^2\right) (\mbox{because $2\leq c_l$ for each $l=1,2,3,4$})\\[6mm]
= &\displaystyle
 \prod r_{ij} \left(\frac{1}{4}\sum_{l=1}^4\left(\frac{d_3(r_{ij},r_{ik},r_{jk})+8r_{ij}r_{ik}r_{jk}}{r_{ij}r_{ik}r_{jk}}\right)^2\right)
\end{array}
\]
and this verifies the weak Four Points Conjecture for cyclic quadrilaterals.\qed

%%%%%%%%%%%%%%%%%%%%%%%%%%%%%%%%%%%%%%%%%%%%%%%%%%%%%%%%%%%%%%%%%%%%%%%%%%%%%%%%%%%%%%%%%%%%%%%%%%%%%%%%%%

%%%%%%%%%%%%%%%%%%%%%%%%%%%%%%%%%%%%%%%%%%%%%%%%%%%%%%%%%%%%%%%%%%%%%%%%%%%%%%% [10]

\section{Almost collinear configurations. {\DJ}okovi\' c's approach}

\subsection{Type (A) configurations}

By a type (A) configurations of $N$ points $x_1,\ldots,x_N$ we shall
mean the case when $N-1$ of the points $x_1,\ldots,x_N$ are
collinear. Set $n=N-1$. In (\cite{DZ}) {\DJ}okovi\' c has proved,
for configurations of type (A), both the Atiyah conjecture
(Theorem 2.1) and the first Atiyah--Sutcliffe conjecture (Theorem
3.1). By using Cartesian coordinates, with $x_i=(a_i,0)$,
$a_1<a_2<\cdots <a_n$ and $x_N=x_{n+1}=(0,b)'$ (with $b=1$), the normalized
Atiyah matrix $M_{n+1}=M_{n+1}(\lambda_1,\ldots,\lambda_n)$
(denoted by $P$ in \cite{DZ} when $b=-1$) is given by
\[
M_{n+1}= \left[
\begin{array}{cccccc}
  1 & \lambda_1 & 0 & \cdots & 0 & 0 \\
  0 & 1 & \lambda_2 & \cdots & 0 & 0 \\
  0 & 0 & 1 & & 0 & 0\\
  \vdots & \vdots & \vdots & \ddots & \vdots & \vdots \\
  0 & 0 &  &  & 1 & \lambda_n \\
  (-1)^{n}e_{n} & (-1)^{n-1}e_{n-1} & \cdots & \cdots & -e_1 & 1
\end{array}
\right]
\]
where $\lambda_1=a_1+\sqrt{a_1^2+b^2}$ $<$
$\lambda_2=a_2+\sqrt{a_2^2+b^2}$ $<\cdots<$
$\lambda_n=a_n+\sqrt{a_n^2+b^2}$ (with $b=1$) are positive real numbers and
where $e_k=e_k(\lambda_1,\ldots,\lambda_n)$, $1\leq k\leq n$, is
the $k$--th elementary symmetric function of
$\lambda_1,\lambda_2,\ldots,\lambda_n$. Its determinant satisfies the inequality
\[
\begin{array}{rl}
  D_n =& 1+\lambda_ne_1+\lambda_n\lambda_{n-1}e_2+\cdots+\lambda_n\lambda_{n-1}\cdots\lambda_1e_n\\
  \geq & 1+e_1(\lambda_1^2,\ldots,\lambda_n^2)+e_2(\lambda_1^2,\ldots,\lambda_n^2)+
  \cdots+e_n(\lambda_1^2,\ldots,\lambda_n^2) \\
  = & \prod_{i=1}^n(1+\lambda_i^2) \\
\end{array}
\]
equivalent to the first Atiyah--Sutcliffe conjecture (\cite{AS},Conjecture 2). The second
Atiyah--Sutcliffe conjecture (\cite{AS},Conjecture 3) for configurations of type (A) is equivalent
to the following inequality
\begin{equation}
\label{E:2.1}
[D_{n+1}(\lambda_1,\ldots,\lambda_n)]^{n-1}\geq
\prod_{k=1}^nD_n(\lambda_1,\ldots,\lambda_{k-1},\lambda_{k+1},\ldots,\lambda_n)
\end{equation}
For $n=2$ this inequality takes the form
\[
1+\lambda_2e_1(\lambda_1,\lambda_2)+\lambda_1\lambda_2e_2(\lambda_1,\lambda_2)
\geq
(1+\lambda_2e_1(\lambda_2))(1+\lambda_1e_1(\lambda_1)
\]
i.e.
\begin{equation}
\label{E:2.2}
1+\lambda_2e_1(\lambda_1,\lambda_2)+\lambda_1\lambda_2e_2(\lambda_1,\lambda_2)
\geq
(1+\lambda_2^2)(1+\lambda_1^2).
\end{equation}
This reduces to $(\lambda_2-\lambda_1)\lambda_1\geq0$, so it is true.
%%%%%%%%%%%%%%%%%%%%%%%%%%%%%%%%%%%%%%%%%%%%%%%%%%%%%%%%%%%%%%%%%%%%%%%%%%%%%%% [11]

Even for $n=3$ the inequality (\ref{E:2.1}) is quite messy thanks
to nonsymmetric character of both sides. Knowing that sometimes it
is easier to solve a more general problem we followed that path
(although we didn't solve the problem in full generality). So let
us start with the case $n=2$. If we look at the following
inequality
\[
1+X_1(\xi_1+\xi_2)+X_1X_2\xi_1\xi_2\geq (1+X_1\xi_1)(1+X_2\xi_2)
\]
which is clearly true if $X_1\geq X_2\geq 0$
and $\xi_1,\xi_2\geq 0$
we obtain the inequality
(\ref{E:2.2}) simply by a specialization $X_1=\xi_1=\lambda_2$, $X_2=\xi_2=\lambda_1$.
So we proceed as follows:

%%%%%%%%%%%%%%%%%%%%%%%%%%%%%%%%%%%%%%%%%%%%%%%%%%%%%%%%%%%%%%%%%%%%%%%%%%%%%%%%%%%%%%%%%%%%%%%%%%%%%%%%%%

Let $\xi_1,\ldots,\xi_n, X_1,\ldots, X_n, n\geq 1$  be two sets of
commuting indeterminates. For any $l, 1\leq l\leq n$ and any
sequences $1\leq i_1\leq\cdots\leq i_l\leq n, 1\leq
j_1,\ldots, j_l\leq n$ we define polynomials
$\Psi^I_J=\Psi^{i_1\ldots i_l}_{j_1\ldots j_l}\in\Q[\xi_1,\ldots,\xi_n,
X_1,\ldots, X_n]$ as follows:
\[
\Psi^I_J:=\sum_{k=0}^le_k(\xi_{j_1},\xi_{j_2},\ldots,\xi_{j_l})X_{i_1}X_{i_2}\cdots
X_{i_k},\ (l\geq 1),\ \Psi^\emptyset_\emptyset:=1\ (j=0)
\]
where $e_k$ is the $k$-th elementary symmetric function.

In particular we have
\[
\begin{array}{ll}
\Psi^i_j & =1+\xi_jX_i,\\[2mm]
\Psi^{i_1i_2}_{j_1j_2}& = 1+(\xi_{j_1}+\xi_{j_2})X_{i_1}+\xi_{j_1}\xi_{j_2}X_{i_1}X_{i_2},\\[2mm]
\Psi^{i_1i_2i_3}_{j_1j_2j_3} & =1+(\xi_{j_1}+\xi_{j_2}+\xi_{j_3})X_{i_1}+
(\xi_{j_1}\xi_{j_2}+\xi_{j_1}\xi_{j_3}+\xi_{j_2}\xi_{j_3})X_{i_1}X_{i_2}+\\
& \hphantom{= 1}+\xi_{j_1}\xi_{j_2}\xi_{j_3}X_{i_1}X_{i_2}X_{i_3},\\
\mbox{ etc.}
\end{array}
\]
The polynomials $\Psi^I_J$ are symmetric w.r.t.\
$\xi_{j_1},\xi_{j_2},\ldots,\xi_{j_l}$, but nonsymmetric w.r.t.\
$X_{i_1},X_{i_2},\ldots,X_{i_l}$.
By specializing $X_{i}$'s to assume real values such that
$X_{i_1}\geq X_{i_2}\geq \ldots\geq X_{i_l}\geq 0$ then we obtain polynomials
in $\xi_j$'s satisfying the following simple but important property.
%%%%%%%%%%%%%%%%%%%%%%%%%%%%%%%%%%%%%%%%%%%%%%%%%%%%%%%%%%%%%%%%%%%%%%%%%%%%%%%%%%
\begin{proposition} (Partition property)\\
\label{R:partition property}
Let $(I_1,\ldots,I_s)$ and $(J_1,\ldots,J_s)$ be ordered set partitions of
respective sets $I=\bigcup_{p=1}^s I_p$ and $J=\bigcup_{p=1}^s J_p$
such that $|I_p|=|J_p|$, $1\leq p\leq s$.
Then the inequality
\[
\Psi^{I}_J\geq\prod_{p=1}^s\Psi^{I_p}_{J_p}
\]
holds coefficientwise w.r.t.\ $\xi_j$'s.
\end{proposition}
\proof
Proof is evident from the definition of $\Psi^I_J$ and the monotonicity of $X_i$'s.\qed

%%%%%%%%%%%%%%%%%%%%%%%%%%%%%%%%%%%%%%%%%%%%%%%%%%%%%%%%%%%%%%%%%%%%%%%%%%%%%%%%%%
For the powers $\left(\Psi^{I}_J\right)^m$ we have the following conjecture.
%%%%%%%%%%%%%%%%%%%%%%%%%%%%%%%%%%%%%%%%%%%%%%%%%%%%%%%%%%%%%%%%%%%%%%%%%%%%%%%%%%
%But there are some more intriguing inequalities with multiset generalization of the
%above partition property.
%%%%%%%%%%%%%%%%%%%%%%%%%%%%%%%%%%%%%%%%%%%%%%%%%%%%%%%%%%%%%%%%%%%%%%%%%%%%%%%%%%
\begin{conjecture}(Weighted Multiset Partition Conjecture)\\
For given natural number $m$ and sets $I$ and $J$, $|I|=|J|$, of natural numbers
let $(I_1,\ldots,I_s)$ and $(J_1,\ldots,J_s)$ be the partitions of the
multiset $I^m$ consisting of $m$ copies of all elements of $I$ and similarly for $J^m$.
\begin{description}
  \item[(i)] Then the inequality
    \[
    \left(\Psi^I_J\right)^m\geq\prod_{p=1}^s\Psi^{I_p}_{J_p}
    \]
    holds coefficientwise w.r.t.\ $\xi_j$'s.
  \item[(ii)] The difference
    \[
    \left(\Psi^I_J\right)^m-\prod_{p=1}^s\Psi^{I_p}_{J_p}
    \]
    is multi--Schur positive with respect to partial alphabets corresponding to the atoms of the
    intersection lattice of the set system $\{J_1,\ldots,J_s\}$.
\end{description}
\end{conjecture}
%%%%%%%%%%%%%%%%%%%%%%%%%%%%%%%%%%%%%%%%%%%%%%%%%%%%%%%%%%%%%%%%%%%%%%%%%%%%%%%%%%
For example, by Partition property, we have the following inequalities
\[
\Psi^{1\ldots n}_{1\ldots n}\geq \Psi^{k}_{k}\Psi^{1..\widehat{k}..n}_{1..\widehat{k}..n},\ (1\leq k\leq n)
\]
which imply the following inequality
\[
\left(\Psi^{1\ldots n}_{1\ldots n}\right)^n
\geq
\prod_{k=1}^{n}\Psi^{k}_{k}\prod_{k=1}^{n}\Psi^{1..\widehat{k}..n}_{1..\widehat{k}..n}
\]
By Partition property we also have the following inequality
\[
\Psi^{1\ldots n}_{1\ldots n}
\geq
\prod_{k=1}^{n}\Psi^{k}_{k}
\]
The last two inequalities suggest the validity of the following inequality
\[
\left(\Psi^{1\ldots n}_{1\ldots n}\right)^{n-1}
\geq
\prod_{k=1}^{n}\Psi^{1..\widehat{k}..n}_{1..\widehat{k}..n}
\]
which is far from obvious (see Conjecture \ref{C:1} below) although
it would be a simple consequence of our Weighted Multiset Partition
Conjecture.

This last conjectural inequality is interesting because it generalizes some special cases of not
yet proven conjectures of Atiyah and Sutcliffe on configurations of
points in three dimensional Euclidean space.

Our conjecture reads as follows:

\begin{conjecture}
\label{C:1}
For any $n\geq 1$, let $X_{1}\geq X_{2}\geq
\ldots\geq X_{n}$ $\geq$ $0$,
$\xi_{1},\xi_{2},\ldots,\xi_{n}\geq 0$, be nonnegative
real numbers. Then we have coefficientwise
(w.r.t.\ $\xi_{1},\xi_{2},\ldots,\xi_{n}$) inequality
\[
\left(
\Psi^{12\cdots n}_{12\cdots n}
\right)^{n-1}
\geq
\prod_{k=1}^{n}\Psi^{12\cdots \hat{k}\cdots n}_{12\cdots \hat{k}\cdots n}
\]
where $12\cdots \hat{k}\cdots n$ denotes the sequence $12\cdots (k-1)(k+1)\cdots n$. The equality obviously holds
true iff $X_1=X_2=\cdots =X_n$.
\end{conjecture}

This Conjecture implies the strongest
Atiyah--Sutcliffe's conjecture for almost collinear configurations
of points (all but one point are collinear, called type(A) in
\cite{DZ}).

To illustrate the Conjecture (\ref{C:1}) we consider first the cases $n=2$ and $n=3$.

\begin{description}
    \item[Case $n=2$:] We have
    \[
        \begin{array}{l@{}l}
          \Psi^{12}_{12} & = 1+(\xi_1+\xi_2)X_1+\xi_1\xi_2X_1X_2=\\[2mm]
                         & = 1+\xi_1X_1+\xi_2X_2+\xi_1\xi_2X_1X_2+(X_1-X_2)\xi_2=\\[2mm]
                         & = (1+\xi_1X_1)(1+\xi_2X_2)+\xi_2(X_1-X_2)\geq\\[2mm]
                         & \geq (1+\xi_1X_1)(1+\xi_2X_2)=\Psi^1_1\Psi^2_2.
        \end{array}
    \]
    \item[Case $n=3$:] We first write $\Psi^{123}_{123}$ in two different ways:
    \[
        \Psi^{123}_{123}=\xi_2(X_1-X_2)+\widehat{\Psi}^{123}_{1\underline{2}3} \mbox{\ \ \  and\ \ \  }
        \Psi^{123}_{123}=\xi_3(X_1-X_2)+\widehat{\Psi}^{123}_{12\underline{3}}.
    \]
    Note that $\widehat{\Psi}^{123}_{1\underline{2}3}$ is obtained from $\Psi^{123}_{123}$ by replacing
    the linear term $\xi_2X_1$ by $\xi_2X_2$, hence all its coefficients are nonnegative.

    The left hand side of the Conjecture (\ref{C:1}) $L_3$ can be rewritten as follows:
    \[
\begin{array}{l@{}l}
  L_3=(\Psi^{123}_{123})^2 & = (\xi_2(X_1-X_2)+\widehat{\Psi}^{123}_{1\underline{2}3})\Psi^{123}_{123}\\[2mm]
                   & = \xi_2(X_1-X_2)\Psi^{123}_{123}+\widehat{\Psi}^{123}_{1\underline{2}3}\Psi^{123}_{123}\\[2mm]
    & = \xi_2(X_1-X_2)\Psi^{123}_{123}+
        \widehat{\Psi}^{123}_{1\underline{2}3}(\xi_3(X_1-X_2)+\widehat{\Psi}^{123}_{12\underline{3}})\\[2mm]
    & = L'_3(X_1-X_2)+\widehat{\Psi}^{123}_{1\underline{2}3}\widehat{\Psi}^{123}_{12\underline{3}}
\end{array}
    \]
    where $L'_3=\xi_2\Psi^{123}_{123}+\xi_3\widehat{\Psi}^{123}_{1\underline{2}3}$ is a positive polynomial.

    Now we have
    \[
    L_3\geq \widehat{L}_3:=\widehat{\Psi}^{123}_{1\underline{2}3}\widehat{\Psi}^{123}_{12\underline{3}}.
    \]
    By using the formula
    \[
    \widehat{\Psi}^{123}_{1\underline{2}3}=\Psi^{12}_{13}+\xi_2X_2\Psi^{13}_{13}=
        (\Psi^2_{2}-1)\Psi^{13}_{13}+\Psi^{12}_{13}
    \]
    we can rewrite $\widehat{L}_3$ as
    \[
\begin{array}{l@{}l}
  \widehat{L}_3 & =\left[ (\Psi^{12}_{13}-\Psi^{13}_{13})+
                           \Psi^2_2\Psi^{13}_{13}\right]\widehat{\Psi}^{123}_{12\underline{3}}\\[2mm]
   & = \xi_1\xi_3X_1(X_2-X_3)\widehat{\Psi}^{123}_{12\underline{3}}+
       \Psi^{13}_{13}(\Psi^2_2\widehat{\Psi}^{123}_{12\underline{3}})
\end{array}
    \]
    The last term in parenthesis can be written as
    \[
\begin{array}{l@{}l}
  \Psi^2_2\widehat{\Psi}^{123}_{12\underline{3}} & = \Psi^{12}_{12}\Psi^{23}_{23}+
                                                     \Psi^1_2(\Psi^{22}_{23}-\Psi^{23}_{23})\\[2mm]
    & = \Psi^{12}_{12}\Psi^{23}_{23}+\xi_2\xi_3X_2(X_2-X_3)\Psi^1_2,
\end{array}
    \]
    so we get
    \[
    \widehat{L}_3=L''_3(X_2-X_3)+\Psi^{12}_{12}\Psi^{13}_{13}\Psi^{23}_{23}
    \]
    where $L''_3$ denotes the positive polynomial
    \[
    L''_3=\xi_1\xi_3X_1\widehat{\Psi}^{123}_{12\underline{3}}+\xi_2\xi_3X_2\Psi^1_2\Psi^{13}_{13}.
    \]
    We now have an explicit formula for $L_3$:
    \[
    L_3=L'_3(X_1-X_2)+L''_3(X_2-X_3)+\Psi^{12}_{12}\Psi^{13}_{13}\Psi^{23}_{23}
    \]
    with $L'_3, L''_3$ positive polynomials, which together with $X_1\geq X_2\geq X_3 (\geq 0)$ implies that
    \[
    L_3\geq R_3:=\Psi^{12}_{12}\Psi^{13}_{13}\Psi^{23}_{23}
    \]
    and the Conjecture (\ref{C:1}) ($n=3$) is proved.
\end{description}

In fact we have proven an instance $n=3$ $\widehat{L}_3\geq R_3$
of a stronger conjecture which we are going to formulate now. Let
$2\leq k\leq n$. We define the modified polynomials
$\widehat{\Psi}^{12\ldots k\ldots n}_{12\ldots \underline{k}\ldots
n}$ as follows:
\[
    \widehat{\Psi}^{12\ldots k\ldots n}_{12\ldots \underline{k}\ldots n}:=
    \xi_k(X_2-X_1)+\Psi^{12\ldots n}_{12\ldots n}
\]
obtained from $\Psi^{12\ldots n}_{12\ldots n}$ by replacing only
one term $\xi_kX_1$ by $\xi_kX_2$, hence $\widehat{\Psi}^{12\ldots
k\ldots n}_{12\ldots \underline{k}\ldots n}$ are still positive.
Let us introduce the following notation:
\[
    \widehat{L}_n:=\prod_{k=2}^{n}\widehat{\Psi}^{12\ldots k\ldots n}_{12\ldots \underline{k}\ldots n}\ ;\ \
    R_n:=\prod_{k=1}^n\Psi^{12\ldots \hat{k}\ldots n}_{12\ldots \hat{k}\ldots n}.
\]
Then clearly $L_n:=(\Psi^{12\ldots n}_{12\ldots n})^{n-1}\geq \widehat{L}_n$.
Now our stronger conjecture reads as
\begin{conjecture}
\label{C:2}
\[\displaystyle \widehat{L}_n\geq R_n\  (n\geq 1)\]
with equality iff $X_2=X_3=\cdots =X_n$.
\end{conjecture}
More generally, we conjecture that the difference $\widehat{L}_n-R_n$ is a polynomial in the
differences $X_2-X_3$, $X_3-X_4$, $\ldots$, $X_{n-1}-X_{n}$ with coefficients in
$\Z_{\geq 0}[X_1,\ldots,X_n,\xi_1,\ldots,\xi_n]$.
\begin{proposition}
\label{P:1}
\[L_n=L'_n(X_1-X_2)+\widehat{L}_n\]
for some positive polynomial $L'_n$.
\end{proposition}
\proof[of Proposition \ref{P:1}]
\[
\begin{array}{l}
 L_n=(\Psi^{12\cdots n}_{12\cdots n})^{n-1}=(\xi_2(X_1-X_2)+\widehat{\Psi}^{12\cdots n}_{1\underline{2}\cdots n})
        (\Psi^{12\cdots n}_{12\cdots n})^{n-2}\\[2mm]
=\xi_2(X_1-X_2)(\Psi^{12\cdots n}_{12\cdots n})^{n-2}+
        \widehat{\Psi}^{12\cdots n}_{1\underline{2}\cdots n}(\xi_3(X_1-X_2)+
        \widehat{\Psi}^{123\cdots n}_{12\underline{3}\cdots n})(\Psi^{12\cdots n}_{12\cdots n})^{n-3}\\[2mm]
=\xi_2(X_1-X_2)(\Psi^{12\cdots n}_{12\cdots n})^{n-2}+
\xi_3(X_1-X_2)\widehat{\Psi}^{12\cdots n}_{1\underline{2}\cdots n}(\Psi^{12\cdots n}_{12\cdots n})^{n-3}+\\
\hphantom{=}+\widehat{\Psi}^{12\cdots n}_{1\underline{2}\cdots n}\widehat{\Psi}^{123\cdots n}_{12\underline{3}\cdots n}
(\Psi^{12\cdots n}_{12\cdots n})^{n-3}\\
\vdots     \\
=(\sum_{k=1}^{n-1}\xi_{k+1}(\prod_{j=2}^k\widehat{\Psi}^{12\ldots j\ldots n}_{12\ldots \underline{j}\ldots n})
(\Psi^{12\ldots n}_{12\ldots n})^{n-k})(X_1-X_2)+
\prod_{j=2}^{n}\widehat{\Psi}^{12\ldots j\ldots n}_{12\ldots \underline{j}\ldots n}.
\end{array}
\]
\qed

Now we turn to study the quotient
\[
\frac{L_n}{R_n}=\frac{(\Psi^{1\ldots n}_{1\ldots n})^{n-1}}
{\displaystyle \prod_{k=1}^{n}\Psi^{1\ldots \widehat{k}\ldots n}_{1\ldots \widehat{k}\ldots n}}
\]
by studying the growth behaviour of quotients of its factors
$\Psi^{1\ldots n}_{1\ldots n}/
{\Psi^{1\ldots \widehat{k}\ldots n}_{1\ldots \widehat{k}\ldots n}}$
w.r.t.\ any of its arguments $X_r$, $1\leq r\leq n$.

%%%%%%%%%%%%%%%%%%%%%%%%%%%%%%%%%%%%%%%%%%%%%%%%%%%%%%%%%%%%%%%%%%%%%%%%%%%%%%%%%%%%%%%%%%%%%%%%%%%%%%%%%%%%%%
%\input{20040127}
\newpage
%\input{20031110}
%%%%%%%%%%%%%%%%%%%%%%%%%%%%%%%%%%%%%%%%%%%%%%%%%%%%%%%%%%%%%%%%%%%%%%%%%%%%%%%%%%%%%%%%%%%%%%%%%%%%%%%%%%%%%%

In the following theorem we obtain an explicit formula for the
numerators of the derivatives w.r.t.\ $X_r, (1\leq r\leq n, r\neq
k)$ of the quantities $\displaystyle {\dogston/\dogstonwok}$. From
this formulas we get some monotonicity properties which enable us to
state some new (refined) conjectures later on.
\begin{theorem}
\label{T:1}
Let
  \begin{equation}
  \label{Eq:0r}
  \Delta_r:=\partial_{X_r}\dogs{1\ldots n}{1\ldots n}\cdot\dogsa{1\ldots\widehat{k}\ldots n}-
            \dogsa{1\ldots n}\cdot\partial_{X_r}\dogsa{1\ldots\widehat{k}\ldots n},
            \ \ (1\leq r\leq n).
  \end{equation}

Then we have the following explicit formulas
\begin{description}
\item[(i)] for any $r$, $1\leq r<k(\leq n)$ we have
\[
    \begin{array}{ll}
    \Delta_r= &\displaystyle \xi_k\!\!\!\!\!\!\!\!\!\!\sum_{0\leq i< r\leq j\leq n}\!\!\!\!\!\!
    s^{(k)}_{(2^i1^{j-i-1})}X_1^2\cdots X_i^2X_{i+1}\cdots \widehat{X}_k\cdots X_j+\\
    &\displaystyle+\!\!\!\!\!\!\!\!\!\!\sum_{0\leq i<r,k\leq j<n}\!\!\!\!\!\!\!\!
    e_ie_j^{(k)}X_1^2\cdots X_i^2X_{i+1}\cdots \widehat{X}_r\cdots \widehat{X}_k\cdots X_j(X_k-X_{j+1})
    \end{array}
  \]
\item[(ii)] for any $r$, $(1\leq )k<r\leq n$ we have
\[
\begin{array}{ll}
\Delta_r=&\displaystyle -\left(
\sum_{0\leq i< r\leq j\leq n}\!\!\!\!\!\!
    s^{(k)}_{(2^i1^{j-i-1})}X_1^2\cdots X_i^2X_{i+1}\cdots \widehat{X}_k\cdots \widehat{X}_r\cdots X_j+\right.\\
    &\displaystyle\left. +\!\!\!\!\!\!\!\!\!\!\sum_{0\leq i<k,r\leq j<n}\!\!\!\!\!\!\!\!
    e_i^{(k)}e_jX_1^2\cdots X_i^2X_{i+1}\cdots \widehat{X}_k\cdots \widehat{X}_r\cdots X_j(X_{j+1}-X_k)
\right)
    \end{array}
\]
where $s^{(k)}_\lambda$ denotes the $\lambda$--th Schur function of
$\xi_1,\ldots ,\xi_{k-1},\xi_{k+1},\ldots,\xi_n$ ($\xi_k$ omitted).
\end{description}

\end{theorem}
\proof[of Theorem \ref{T:1}]

  {\bf (i)} For any $r$, $1\leq r<k(\leq n)$ we find explicitly a formula
  as follows. We shall use notations $X_\twodots{i}:=X_1X_2\cdots X_i$,
  for multilinear monomials and $e_i:=e_i(\xi_1,\ldots,\xi_n)$,
  $e_i^{(k)}=e_i(\xi_1,\ldots,\widehat{\xi_k},\ldots\xi_n)$ for the elementary
  symmetric functions (here $k$ is fixed). Then we can rewrite our basic quantities
  \begin{equation}
  \label{Eq:1r}
  \dogsa{1\ldots n}:=\sum_{i=0}^{n}e_iX_\twodots{i}
  \end{equation}
  \begin{equation}
  \label{Eq:2r}
  \begin{array}{lrl}
  \dogsa{1\ldots\widehat{k}\ldots n}&:=&\displaystyle\sum_{i=0}^{k-1}e_i^{(k)}X_\twodots{i}+
  \frac{1}{X_k}\sum_{i=k}^{n-1}e_i^{(k)}X_\twodots{i+1}=\\[4mm]
    &=&\displaystyle\sum_{i=0}^{n-1}e_i^{(k)}X_\twodots{i}+
    \frac{1}{X_k}\sum_{i=k}^{n-1}e_i^{(k)}X_\twodots{i}(X_{i+1}-X_k)
  \end{array}
  \end{equation}
  For the derivatives we get immediately
  \begin{equation}
  \label{Eq:3r}
  \partial_{X_r}\dogston=\frac{1}{X_r}\sum_{i=r}^ne_iX_\twodots{i}=\frac{1}{X_r}
  \left(\dogston-\sum_{i=0}^{r-1}e_iX_\twodots{i}\right)
  \end{equation}
  \begin{equation}
  \label{Eq:4r}
    \partial_{X_r}\dogstonwok=\frac{1}{X_r}\sum_{i=r}^{n-1}e_i^{(k)}X_\twodots{i}+
    \frac{1}{X_kX_r}\sum_{i=k}^{n-1}e_i^{(k)}X_\twodots{i}(X_{i+1}-X_k)
  \end{equation}
  \begin{equation}
  \label{Eq:5r}
    \hphantom{\partial_{X_r}\dogstonwok}=\frac{1}{X_r}\left(\dogstonwok-\sum_{i=0}^{r-1}e_i^{(k)}X_\twodots{i}\right)
  \end{equation}
  By plugging (\ref{Eq:3r}) and (\ref{Eq:5r}) into (\ref{Eq:0r}) we obtain
\[
X_r\Delta_r =\dogston\left( \sum_{i=0}^{r-1}e_i^{(k)}X_\twodots{i}\right)-
                            \dogstonwok\left( \sum_{i=0}^{r-1}e_iX_\twodots{i}\right)=
\]
and after simple cancelation, by invoking (\ref{Eq:2r}) we get
\[
\begin{array}{l}
=\left(\sum_{j=r}^{n}e_jX_\twodots{j}\right)\left(\sum_{i=0}^{r-1}e_i^{(k)}X_\twodots{i}\right)-\\
\hspace{5ex}
\left(\sum_{j=r}^{n-1}e_j^{(k)}X_\twodots{j}+
\frac{1}{X_k}\sum_{j=k}^{n-1}e_j^{(k)}X_\twodots{j}(X_{j+1}-X_k)\right)
\left(\sum_{i=0}^{r-1}e_iX_\twodots{i}\right)
\end{array}
\]
i.e.
{\small
  \[
    X_r\Delta_r=\!\!\!\!\!\!\sum_{0\leq i< r\leq j\leq n}\!\!\!\!\!\!(e_je_i^{(k)}-
    e_ie_j^{(k)})X_\twodots{i}X_\twodots{j}
    +\frac{1}{X_k}\sum_{0\leq i<r,k\leq j<n}\!\!\!\!\!\!\!\!e_ie_j^{(k)}X_\twodots{i}X_\twodots{j}(X_k-X_{j+1})
  \]
}
  If we use a simple identity $e_j=e_j^{(k)}+\xi_ke_{j-1}^{(k)}$, we can identify the quantity
  \[
  \begin{array}{l}
    e_je_i^{(k)}-e_ie_j^{(k)}=(e_j^{(k)}+\xi_ke_{j-1}^{(k)})e_i^{(k)}-(e_i^{(k)}+\xi_ke_{i-1}^{(k)})e_j^{(k)}=\\[2ex]
    \hphantom{e_je_i^{(k)}-e_ie_j^{(k)}}
    =\left|
    \begin{array}{cc}
    e_{j-1}^{(k)}&e_j^{(k)}\\
    e_{i-1}^{(k)}&e_i^{(k)}
    \end{array}
    \right|\xi_k
    =s^{(k)}_{2^i1^{j-i-1}}\xi_k
  \end{array}
  \]
  Thus in this case $(1\leq r<k)$ we obtain a formula
  \[
    \begin{array}{ll}
    \Delta_r= &\displaystyle \xi_k\!\!\!\!\!\!\!\!\!\!\sum_{0\leq i< r\leq j\leq n}\!\!\!\!\!\!
    s^{(k)}_{(j-1,i)}X_1^2\cdots X_i^2X_{i+1}\cdots \widehat{X}_k\cdots X_j+\\
    &\displaystyle\hphantom{xxxx}+\!\!\!\!\!\!\!\!\!\!\sum_{0\leq i<r,k\leq j<n}\!\!\!\!\!\!\!\!
    e_ie_j^{(k)}X_1^2\cdots X_i^2X_{i+1}\cdots \widehat{X}_r\cdots \widehat{X}_k\cdots X_j(X_k-X_{j+1})
    \end{array}
  \]
  (where $e_j^{(k)}=e_j^{(k)}=e_j(\xi_1,\ldots ,\widehat{\xi_k},\ldots,\xi_n)$)
  in terms of Schur functions (of arguments $\xi_1,\ldots ,\widehat{\xi_k},\ldots,\xi_n$)
  corresponding to a transpose $(2^i1^{j-i-1})$ of a partition $(j-1,i)$
  (cf.\ Jacobi--Trudi formula, I 3.5 in \cite{Jacobi}).\\[5mm]

%%%%%%%%%%%%%%%%%%%%%%%%%%%%%%%%%%%%%%%%%%%%%%%%%%%%%%%%%%%%%%%%%%%%%%

  {\bf (ii)} For any $r$, $(1\leq )k<r\leq n$. In this case we use
  \[
\partial_{X_r}\dogstonwok=\frac{1}{X_kX_r}\sum_{j=r-1}^{n-1}e_j^{(k)}X_\twodots{j+1}
  \]
\[
\begin{array}{l@{\ =\ }l}
\dogstonwok
&\displaystyle
\sum_{i=0}^{k-1}e_i^{(k)}X_\twodots{i}+\frac{1}{X_k}\sum_{i=k}^{n-1}e_i^{(k)}X_\twodots{i+1}=\\
&\displaystyle
\frac{1}{X_k}\left(\sum_{i=0}^{k-1}X_\twodots{i}(X_k-X_{i+1})+\sum_{i=0}^{n-1}e_i^{(k)}X_\twodots{i}\right)
\end{array}
\]
By plugging this into (\ref{Eq:0r}) we get
\[
\begin{array}{ll}
X_kX_r\Delta_r
&\displaystyle =\left( \sum_{j=r}^{n}e_jX_\twodots{j}\right)
 \left( \sum_{i=0}^{k-1}e_i^{(k)}X_\twodots{i}(X_k-X_{i+1})+\sum_{i=0}^{n-1}e_i^{(k)}X_\twodots{i+1}\right)-\\
 &\displaystyle \hphantom{xxxx}-\left(\sum_{j=0}^{r-1}e_jX_\twodots{j}+\sum_{j=r}^{n}e_jX_\twodots{j}\right)
 \left(\sum_{i=r-1}^{n-1}e_i^{(k)}X_\twodots{i+1}\right)\\
&\hspace{-10ex}\displaystyle
 =\left( \sum_{i=0}^{r-2}e_i^{(k)}X_\twodots{i+1}\right)\left( \sum_{j=r}^{n}e_jX_\twodots{j}\right)-
\left( \sum_{i=0}^{r-1}e_iX_\twodots{i}\right)\left(\sum_{j=r-1}^{n-1}e_j^{(k)}X_\twodots{j+1} \right)+\\
 &\displaystyle \hphantom{xxxx}+\sum_{i=0}^{k-1}\sum_{j=r}^{n}e_i^{(k)}e_jX_\twodots{i}X_\twodots{j}(X_k-X_{i+1})\\
&\hspace{-10ex}\displaystyle
=\left( \sum_{i=1}^{r-1}e_{i-1}^{(k)}X_\twodots{i}\right)\left( \sum_{j=r}^{n}e_jX_\twodots{j}\right)-
\left( \sum_{i=0}^{r-1}e_iX_\twodots{i}\right)\left(\sum_{j=r}^{n}e_{j-1}^{(k)}X_\twodots{j} \right)+\\
 &\displaystyle \hphantom{xxxx}+\sum_{i=0}^{k-1}\sum_{j=r}^{n}e_i^{(k)}e_jX_\twodots{i}X_\twodots{j}(X_k-X_{i+1})\\
\end{array}
\]
By using a formula for elementary symmetric functions ($e_i=e_i^{(k)}+\xi_ke_{i-1}^{(k)}$)
we can write in terms of Schur functions (of arguments
$\xi_1,\ldots ,\xi_{k-1},\xi_{k+1},\ldots,\xi_n$), where ${\lambda '}$ is a conjugate
of $\lambda$.
\[
e_{i-1}^{(k)}e_j-e_ie_{j-1}^{(k)}=e_{i-1}^{(k)}e_j^{(k)}-e_i^{(k)}e_{j-1}^{(k)}=
-\left|
\begin{array}{cc}
e_{j-1}^{(k)} & e_j^{(k)}\\
e_{i-1}^{(k)} & e_i^{(k)}
\end{array}
\right|=-s^{(k)}_{2^i1^{j-i-1}}=-s^{(k)}_{(j-1,i)'}
\]
Thus we obtain a formula
\[
\begin{array}{ll}
\Delta_r=&\displaystyle -\left(
\sum_{0\leq i< r\leq j\leq n}\!\!\!\!\!\!
    s^{(k)}_{(j-1,i)'}X_1^2\cdots X_i^2X_{i+1}\cdots \widehat{X}_k\cdots \widehat{X}_r\cdots X_j+\right.\\
    &\displaystyle\hphantom{xxxx}\left. +\!\!\!\!\!\!\!\!\!\!\sum_{0\leq i<k,r\leq j<n}\!\!\!\!\!\!\!\!
    e_i^{(k)}e_jX_1^2\cdots X_i^2X_{i+1}\cdots \widehat{X}_k\cdots \widehat{X}_r\cdots X_j(X_{j+1}-X_k)
\right)
    \end{array}
\]
\begin{corollary} ($X_r$--monotonicity)\\
\label{C:2.1}
Let $X_1\geq\cdots\geq X_n\geq 0$, $\xi_1,\ldots,\xi_n\geq 0$ be as before. Then
\begin{enumerate}[(i)]
  \item for any $r$, $1\leq r<k\ (\leq n)$ we have
  \[
\frac{\Psi^{1\ldots n}_{1\ldots n}}{
{\Psi^{1\ldots \widehat{k}\ldots n}_{1\ldots \widehat{k}\ldots n}}}
\geq
\frac{\Psi^{1\ldots\ r+1\ r+1\ \ldots n}_{1\ldots\ \hspace{1.2ex}r\hspace{2ex}r+1\ \ldots n}}{
{\Psi^{1\ldots\ r+1\ r+1\ \ldots\widehat{k}\ldots n}_{1\ldots\ \hspace{1.2ex}r\hspace{2ex}r+1\ \ldots \widehat{k}\ldots n}}}
  \]

  \item for any $r$, $(1\leq)\ k<r\ (\leq n)$ we have
  \[
\frac{\Psi^{1\ldots n}_{1\ldots n}}{
{\Psi^{1\ldots \widehat{k}\ldots n}_{1\ldots \widehat{k}\ldots n}}}
\geq
\frac{\Psi^{1\ldots\ r-1\ r-1\ \ldots n}_{1\ldots\ r-1\hspace{2ex}r\ \hspace{1.2ex}\ldots n}}
{{\Psi^{1\ldots\widehat{k}\ldots\ r-1\ r-1\ \ldots n}_{1\ldots\widehat{k}\ldots\ r-1\hspace{2ex}r\ \hspace{1.2ex}\ldots n}}}
  \]
\end{enumerate}
\end{corollary}
Now we illustrate how to use Corollary \ref{C:2.1} to prove our Conjecture \ref{C:1} for $n=2, 3, 4$ and $5$.\\
{\bf Case $n=2$}
\[
Q_2:=\frac{\dogsa{12}}{\dogsa{1}\dogsa{2}}\geq\frac{\dogs{22}{12}}{\dogs{2}{1}\dogsa{2}}=1\mbox{ (by $(i)$)}
\]
{\bf Case $n=3$}
\[
\begin{array}{rl}
\displaystyle
Q_3:=\frac{\dogsa{123}\dogsa{123}}{\dogsa{12}\dogsa{13}\dogsa{23}}&\displaystyle
\geq
\frac{\dogs{223}{123}\dogsa{123}}{\dogs{22}{12}\dogsa{13}\dogsa{23}}
\geq
\frac{\dogs{223}{123}\dogs{223}{123}}{\dogs{22}{12}\dogsa{13}\dogsa{23}}\mbox{ (by $(i)$)}\\[5mm]
&\displaystyle
\geq
\frac{\dogs{222}{123}\dogs{223}{123}}{\dogs{22}{12}\dogs{22}{13}\dogsa{23}}
\geq
\frac{\dogs{222}{123}\dogs{222}{123}}{\dogs{22}{12}\dogs{22}{13}\dogsa{23}}=1\mbox{ (by $(ii)$)}
\end{array}
\]
{\bf Case $n=4$}
\[
\begin{array}{rl}
\displaystyle
Q_4:=\frac{(\dogsa{1234})^3}{\dogsa{123}\dogsa{124}\dogsa{134}\dogsa{234}}
\geq\cdots\geq
\frac{\dogs{2244}{1234}(\dogs{2224}{1234})^2}{\dogs{224}{123}\dogs{224}{124}\dogs{224}{134}\dogs{224}{234}}
\ \ (\geq 1)
\end{array}
\]
This last inequality follows from the following symmetric function identity:
\[
\begin{array}{l}
\dogs{2244}{1234}(\dogs{2224}{1234})^2-\dogs{224}{123}\dogs{224}{124}\dogs{224}{134}\dogs{224}{234}=\\[2mm]
X_2^2X_4^4m_{2222}+2X_2^2X_4^3m_{2221}+X_2^2X_4^2m_{222}+3X_2^2X_4^2m_{2211}+X_2^2X_4m_{221}\\[2mm]
+4X_2^2X_4m_{2111}+X_2^2m_{211}+X_2(3X_2+2X_4)m_{1111}+X_2m_{111}
\end{array}
\]
where $m_\lambda=m_\lambda(\xi_1,\xi_2,\xi_3,\xi_4)$ are the monomial symmetric functions.\\
{\bf Case $n=5$}
\[
\begin{array}{rl}
\displaystyle
Q_5:=\frac{(\dogstona{5})^4}{\prod_{k=1}^5\dogstonwoka{5}}
\geq\cdots\geq
\frac{(\dogs{22244}{12345}\dogs{22444}{12345})^2}
{\dogs{2244}{1234}\dogs{2244}{1235}\dogs{2244}{1245}\dogs{2244}{1345}\dogs{2244}{2345}}
\ \ (\geq 1)
\end{array}
\]
The last inequality is equivalent to an explicit symmetric function identity with all
coefficients (w.r.t.\ monomial basis) positive.

Now we state our stronger conjecture.
\begin{conjecture}(for symmetric functions)\\
\label{C:3.6}
Let $X_1\geq X_2\geq \cdots\geq X_n\geq 0$ and
$\xi_1,\ldots,\xi_n\geq 0$.
Then the inequalities
\begin{enumerate}[(a)]
\item For $n$ even
\[
\!\!\!\!\!\!\!\!\!\!\!\!
\dogs{2\ 2\ 4\ 4\ldots n\ n}{1\ 2\hspace{0.61ex}\ldots\hspace{0.61ex} n-1\ n}
\left(
\prod_{k=1}^{n/2}\dogs{2\ 2\ 4\ 4\ldots 2k\ 2k\ 2k\ldots n-2\ n-2\ n}
{1\ 2\ 3\ 4 \hspace{6.5ex}\ldots\hspace{7.2ex} n-1\ n}
\right)^2
\geq
\prod_{k=1}^n\dogs{2\ 2\ 4\ 4\ldots n-2\ n-2\ n}
{1\ 2 \hspace{1.31ex}\ldots\hspace{1.31ex} \widehat{k}\hspace{1.31ex}\ldots\hspace{1.31ex} n-1\ n}
\]

\item For $n$ odd
\[
\left(
\prod_{k=1}^{\lfloor n/2\rfloor}\dogs{2\ 2\ 4\ 4\ldots 2k\ 2k\ 2k\ldots n-1\ n-1}
{1\ 2\ 3\ 4\hspace{2.35ex}\ldots\hspace{6.8ex} n-1\hspace{2ex} n}
\right)^2
\geq
\prod_{k=1}^n\dogs{2\ 2\ 4\ 4\ldots n-1\ n-1}
{1\ 2 \hspace{1.57ex}\ldots\hspace{1.57ex} \widehat{k}\hspace{1.57ex}\ldots\hspace{1.57ex} n}
\]
\end{enumerate}
hold true coefficientwise ($m$--positivity).
\end{conjecture}

Now we motivate another inequalities for symmetric functions which also refine the strongest
Atiyah--Sutcliffe conjecture for configurations of type (A). Let $n=3$.
We apply Corollary \ref{C:2.1} by using steps $(ii)$ only.
\[
Q_3:=\frac{\dogsa{123}\dogsa{123}}{\dogsa{12}\dogsa{13}\dogsa{23}}
\geq
\frac{\dogs{113}{123}\dogsa{123}}{\dogsa{12}\dogsa{13}\dogs{13}{23}}
\geq
\frac{\dogs{112}{123}\dogsa{123}}{\dogsa{12}\dogs{12}{13}\dogs{13}{23}}
\geq
\frac{\dogs{112}{123}\dogs{122}{123}}{\dogsa{12}\dogs{12}{13}\dogs{12}{23}}
\geq
1
\]
The last inequality is equivalent to nonnegativity of the expression
\[
\dogs{112}{123}\dogs{122}{123}-\dogsa{12}\dogs{12}{13}\dogs{12}{23}\ \
(=X_1(X_1-X_2)^2\xi_1\xi_2\xi_3\geq 0).
\]
Similarly, for $n=4$, the symmetric function inequality stronger than $Q_4\geq 1$
would be the following
\[
\dogs{1123}{1234}\dogs{1223}{1234}\dogs{1233}{1234}
\geq
\dogsa{123}\dogs{123}{124}\dogs{123}{134}\dogs{123}{234}
\]
Now we state a general conjecture for symmetric functions which imply the strongest
Atiyah--Sutcliffe conjecture for almost collinear type (A) configurations.
\begin{conjecture}
\label{C:2.2}
Let $X_1\geq\cdots\geq X_n\geq$, $\xi_1,\ldots\xi_n\geq 0$. Then the following inequality
for symmetric functions in $\xi_1,\ldots,\xi_n$
\[
\dogs{112\ldots n-1}{123\ldots n}
\dogs{1223\ldots n-1}{1234\ldots n}
\cdots
\dogs{12\ldots n-2\ n-1\ n-1}{12\ldots n-2\ n-1\ n}
\geq
\dogsa{1\ 2\ldots n-1}
\dogs{1\ 2\ldots n-1}{1\ 2\ldots n-2\ n}
\cdots
\dogs{1\ 2\ldots n-1}{2\ 3\ldots n-1}
\]
i.e.
\[
\prod_{k=1}^{n-1}
\dogs{1\ 2\ldots k\ \hphantom{+}k\hphantom{1}\ldots n}{1\ 2\ldots k\ k+1\ldots n}
\geq
\prod_{k=1}^{n}
\dogs{1\ 2\ \ldots\  n-1}{1\ 2\ldots \widehat{k}\ldots n}
\]
holds true coefficientwise ($m$--positivity).
\end{conjecture}
\begin{remark}
Conjectures \ref{C:3.6} and \ref{C:2.2} seems to hold also for the Schur
basis of symmetric functions in $\xi_1,\ldots,\xi_n$.
\end{remark}

We have checked this Conjecture \ref{C:2.2} up to $n=5$ by using {\tt Maple} and symmetric
function package {\tt SF} of J.\ Stembridge. For $n$ bigger than five the computations are extremely
intensive and hopefully in the near future would be possible by using more powerful
computers.

Note that the right hand side of the Conjecture \ref{C:2.2} involves symmetric functions
of partial alphabets $\xi_1, \xi_2,\ldots,\xi_{k-1},\xi_{k+1},\ldots,\xi_n$. But the left hand side
doesn't have this "defect". Our objective now is to give explicit formula for the right hand
side in terms of the elementary symmetric functions of the full alphabet
$\xi_1, \xi_2,\ldots,\xi_n$. This we are going to achieve by using resultants as follows.
\begin{lemma}
\label{L:res}
For any $k$, ($1\leq k\leq n$), we have
\[
\dogs{1\ldots k\ldots n-1}{1\ldots \widehat{k}\ldots n}=
\sum_{j=0}^{n-1}a_j\xi_k^{n-1-j}
\]
where
\[
\begin{array}{l}
a_{n-1} = 1+X_1e_1+X_1X_2e_2+\ldots +X_1\cdots X_{n-1}e_{n-1},\\
a_{n-2} = -X_1-X_1X_2e_1-\ldots -X_1\cdots X_{n-1}e_{n-2},\\
\cdots\\
a_0 = (-1)^{n-1}X_1\cdots X_{n-1}\\[5mm]
\mbox{i.e.}\\
\displaystyle a_{n-1-j}=(-1)^j\sum_{i=j}^{n-1}X_1\cdots X_ie_{i-j}
\end{array}
\]
\end{lemma}
\proof[of Lemma \ref{L:res}]
By definition we have
\begin{equation}
\label{E:2.8}
\dogs{1\ldots n-1}{1\ldots\widehat{k}\ldots n}=\sum_{i=0}^{n-1}X_1\cdots X_ie_i^{(k)}
\end{equation}
where $e_i^{(k)}$ is the $i$--th elementary function of $\xi_1,\ldots ,\xi_{k-1},\xi_{k+1},\ldots,\xi_n$.
Now from the decomposition
\[
(1+\xi_kt)^{-1}\prod_{j=1}^{n}(1+\xi_jt)=\prod_{j\neq k}(1+\xi_jt)=
\sum_{i=0}^{n-1}e_i^{(k)}t^i
\]
we get
\[
e_i^{(k)}=e_i-e_{i-1}\xi_k+e_{i-2}\xi_k^2-\cdots +(-1)^{i}\xi_k^i
\]
By substituting this into equation (\ref{E:2.8}) the Lemma \ref{L:res} follows.
\qed\\

Then, by Lemma \ref{L:res}, the right hand side
\[
R_n=\prod_{k=1}^{n}\dogs{1\ 2\ \ldots\ k\ \ldots \  n-1}{1\ 2\ \ldots\  \widehat{k}\ \ldots\ \ \ n}
=
\prod_{k=1}^n\left(\sum_{j=0}^{n-1}a_j\xi_k^{n-1-j} \right)
\]
can be understood as a resultant $R_n=Res(f,g)$ of the following two polynomials
\[
\begin{array}{r@{\ =\ }l}
f(x) & \displaystyle \sum_{j=0}^{n-1}a_jx^{n-1-j}\\[2mm]
g(x) & \displaystyle \prod_{i=1}^{n}(x-\xi_i)=\sum_{j=0}^{n}(-1)^je_jx^{n-j}
\end{array}
\]
The Sylvester formula
\[
R_n=
\left|
\begin{array}{cccccccc}
  1 & -e_1 & e_2 & -e_3 & \ldots & (-1)^{n}e_n &   &  \\
    & 1 & -e_1 & e_2 & -e_3 & \ldots &   &  \\
    &   & \ddots &   &   &   &   &   \\
    &   &   & 1 &  -e_1 & \cdots  &   &   \\
  a_0 & a_1 & a_2 & \cdots & a_n &   &   &   \\
  & a_0 & a_1 & a_2 & \cdots & a_n &   &  \\
    &   & \ddots &   &   &   &   &   \\
    &   &   & a_0 & a_1 & a_2 & \cdots & a_n \\
\end{array}
\right|
\ \
\left(=:
\left|
\begin{array}{cc}
A & B\\
C & D
\end{array}
\right|
\right)
\]
can be simplified as
\[
=|A|\cdot |D-CA^{-1}B|=|D-CA^{-1}B|.
\]
The entries of the $n\times n$ matrix $\Delta:=D-CA^{-1}B$ are given by
\[
\delta_{ij}=\left\{
\begin{array}{ll}
\displaystyle (-1)^{j-i-1}\sum_{k=j+1}^{n}X_1\cdots X_{k+i-j}e_{k}, & 0\leq i<j\leq n-1\\
\displaystyle (-1)^{j-i}\sum_{k=0}^{j}X_1\cdots X_{k+i-j}e_{k}, & 0\leq j\leq i\leq n-1
\end{array}
\right.
\]
%\[
%\delta_{ij}=\left\{
%\begin{array}{ll}
%\displaystyle \sum_{k=0}^{i-j}(-1)^{i+j}X_1\cdots X_{i-k}e_{j-k} & \mbox{ for } n-1\geq i\geq j\geq 0\\
%\displaystyle \sum_{k=0}^{n}(-1)^{j-i}X_1\cdots X_{j-i+k}e_{j-i+k+1} & \mbox{ for } n-1\geq j> i\geq 0
%\end{array}
%\right.
%\]
For example, for $n=3$
\[
\Delta_3=\left|
\begin{array}{ccc}
  1 & X_1e_2+X_1X_2e_3 & -X_1e_3 \\[2mm]
  -X_1 & 1+X_1e_1 & X_1X_2e_3 \\[2mm]
  X_1X_2 & -X_1-X_1X_2e_1 & 1+X_1e_1+X_1X_2e_2
\end{array}
\right|
\]
By elementary operations we get
\[
\Delta_3=
\left|
\begin{array}{ccc}
  1 & * & * \\[2mm]
  0 & \dogs{112}{123} & X_1(X_2-X_1)e_3 \\[2mm]
  0 & X_2-X_1 & \dogs{122}{123}
\end{array}
\right|
=
\left|
\begin{array}{cc}
\dogs{112}{123} & X_1(X_2-X_1)e_3 \\[2mm]
X_2-X_1 & \dogs{122}{123}
\end{array}
\right|
\]
Similarly, for $n=4$ we obtain
{\small
\[
\Delta_4=\left|
\begin{array}{ccc}
  \dogs{1123}{1234} & -X_1(X_1-X_2)e_3-X_1X_2(X_1-X_3)e_4 & X_1(X_1-X_2)e_4 \\[2mm]
  -(X_1-X_2) & \dogs{1223}{1234} & -X_1X_2(X_2-X_3)e_4 \\[2mm]
  X_1(X_2-X_3) & -(X_1-X_3)-X_1(X_2-X_3)e_1 & \dogs{1233}{1234}
\end{array}
\right|
\]
}
In general
\[
\Delta_n=\det(\delta'_{ij})_{1\leq i,j\leq n-1}
\]
where
\[
\delta'_{ij}=\left\{
\begin{array}{l@{\ ,\ }l}
\displaystyle (-1)^{j-i}\sum_{k=j+1}^{n}X_1\cdots X_{k+i-j-1}(X_i-X_{k+i-j})e_{k} & 1\leq i<j\leq n-1\\[3mm]
\displaystyle \dogs{1\ \ldots\ i\ i\ \ldots\ n}{1\ 2\ \ldots\ n} & i=j\\[2mm]
\displaystyle (-1)^{j-i}\sum_{k=0}^{j}X_1\cdots X_{k+i-j-1}(X_{k+i-j}-X_i)e_{k} & 1\leq j<i\leq n-1
\end{array}
\right.
\]
%\[
%\delta'_{ij}=\left\{
%\begin{array}{l@{\ ,\ }l}
%\displaystyle \sum_{k=j+1}^{n-1}(-1)^{i+j+1}X_1\cdots X_{k-1}(X_k-X_i)e_{k+1} & \mbox{ for } i<j\\[2mm]
%\displaystyle \dogs{1\ \ldots\ i\ i\ \ldots\ n-1}{1\ 2\ \ldots\ n} &  \mbox{ for } i=j\\[2mm]
%\displaystyle \sum_{k=0}^{j-1}(-1)^{i+j}X_1\cdots X_{i-k-2}(X_{i-k-1}-X_i)e_{k} & \mbox{ for } i>j
%\end{array}
%\right.
%\]
\begin{corollary}
The conjecture \ref{C:2.2} is equivalent to a Hadamard type
inequality, holding coefficientwise, for the (non Hermitian) matrix $\displaystyle
(\delta'_{ij})_{1\leq i,j\leq n-1}$, i.e.
\[
\prod_{i=1}^{n-1}\delta'_{ii}\geq \det(\delta'_{ij})
\]
\end{corollary}

%%%%%%%%%%%%%%%%%%%%%%%%%%%%%%%%%%%%%%%%%%%%%%%%%%%%%%%%%%%%%%%%%%%%%%%%%%%%%%%%%%%%%%%%%%%%%%%%%%%%%%%%%%

\section{Verification of the {\DJ}okovi\' c's strengthening of the Atiyah--Sutcliffe Conjecture (C2)
for some nonplanar configurations with dihedral symmetry}

Here we basically follow {\DJ}okovi\' c's \cite{DZb}, where only Atiyah
conjecture C1 was proved, make some additional refinements
including a proof of Atiyah--Sutcliffe conjecture C2.

Let $N=m+n$ points be such that
\begin{enumerate}
  \item The first $m$ points $x_1,\ldots,x_m$ lie on a line $L$.
  \item The remaining $n$ points $y_j=x_{m+j+1}$ $(j=0,1,...,n-1)$
  are the vertices oif a regular $n$--gon whose plane is
  perpendicular to $L$ and whose centroid lies on $L$.
\end{enumerate}
We may assume $L=\R\times \{0\}\subset \R\times\C=\R^3$ and write
$x_i=(a_i,0)$, $1\leq i\leq m$, $a_1\leq\ldots\leq a_m$ and
$y_j=(0,b_j)$, $b_j=-\xi^j$, $\xi=e^{2\pi i/n}$, $0\leq j\leq n-1$.

We set
\[
\lambda_i=a_i+\sqrt{1+a_i^2}
\]
Recall that $a_1<\cdots<a_m$ and, consequently
$0<\lambda_1<\cdots<\lambda_m$. Then the associated polynomials
$p_i$ (up to scalar factors) are given by
\[
\begin{array}{l}
\displaystyle p_i(x,y)=x^{m-i}y^{i-1}(x^n-\lambda_i^{n}y^n),\ 1\leq i\leq m\\
\displaystyle p_{m+j+1}(x,y)=\prod_{s\neq
j}\left(x+\frac{\overline{b_s}-\overline{b_j}}{|b_s-b_j|}y\right)
\cdot \prod_{i=1}^m (y-\lambda_ib_jx), 0\leq j<n
\end{array}
\]
By noting that
\[
b_s-b_j=2i\xi^{\frac{j+s}{2}}\sin \frac{\pi (j-s)}{n}
\]
(in {\DJ}okovi\' c $\xi^{j+s}$ should be replaced by
$\xi^{\frac{j+s}{2}}$) we obtain
\[
x+\frac{\overline{b_s}-\overline{b_j}}{|b_s-b_j|}y=
\left(
-\overline{b_j}y-i\xi^{\frac{s-j}{2}}\mbox{sgn}(s-j)
\right)
\frac{1-\overline{b_s}b_j}{|b_s-b_j|}
\]
and
\[
y-\lambda_ib_jx=-b_j(-\overline{b_j}y+\lambda_ix)
\]
Note also that
\[
\{
\xi^\frac{s-j}{2}\mbox{sgn}(s-j)|s=1,\ldots,j-1,j+1,\ldots,n
\}
=
\{
e^{\pi i k/n}|k=1,\ldots,n-1
\}
\]
Thus, after dehomogenizing the polynomials $p_i$ by setting $x=1$, we
obtain (up to scalar factors) the following polynomials:
\[
\begin{array}{l}
\widetilde{P}_i(y)=y^{i-1}(1-\lambda_i^ny^n),\ 1\leq i\leq m;\\
\widetilde{P}_{m+j+1}(y)=f(\xi^{-1}y),\ 0\leq j< n
\end{array}
\]
where
\[
f(y)=\prod_{s=1}^{n-1}(y-ie^{\pi i s/n})\prod_{i=1}^m(y+\lambda_i)
\]
(in {\DJ}okovi\' c the last $n$ polynomials are reordered)

The main result of {\DJ}okovi\' c is the Theorem {\it 3.1} where he proved
Atiyah conjecture for configurations described above, by explicitly
computing the determinant of the coefficients matrix $\widetilde{P}$
of the polynomials $\{\widetilde{p}_k(y)|k=1,\ldots,\underbrace{m+n}_N\}$ in terms
of the coefficients of
\[
f(y)=\sum_{k=0}^{N-1}\widetilde{E}_ky^{N-1-k}
\]
His formula reads as follows:
\[
\left|
\det(\widetilde{P})
\right|
=
n^{n/2}\prod_{k=0}^{n-1}f_k
\]
where
\[
f_k=\sum_{s\geq 0}
\left(
\prod_{j=1}^s\lambda^n_{N-jn-k}
\right)
\widetilde{E}_{k+sn},\ 0\leq k<n.
\]
We shall now present an amazingly simple formula for coefficients of the polynomial
\[
h(y):=\prod_{s=1}^{n-1}(y-ie^{\pi i s/n})=\sum_{j=0}^{n-1}c_jy^{n-1-j}
\]
\begin{proposition}
\label{P:2}
let $\gamma_k:=\cot \left(\frac{k\pi}{2n}\right)$. Then
\[
c_0=1,\ c_j=\prod_{k=1}^{j}\gamma_k\ \ \ \ (1\leq j\leq n-1)
\]
\end{proposition}
\proof
Put $\xi_k=-ie^{\pi i k/n}$, $k=1,\ldots,n-1$. Then
\[
\begin{array}{l}
c_j=\mbox{the j--th elementary symmetric function of }\xi_1,\ldots,\xi_{n-1}\\
\hphantom{c_j}=e_j(\xi_1,\ldots,\xi_{n-1})
\end{array}
\]
Let us first compute the power sums
\[
\begin{array}{l@{\ =\ }l}
\displaystyle
p_s=\sum_{k=1}^{n-1}\xi_k^s=(-i)^s\sum_{k=1}^{n-1}e^{\pi i s k/n}
& (-i)^s(e^{\pi i s /n}-e^{\pi i s})/(1-e^{\pi i s})\\[3mm]
&
\left\{
\begin{array}{l}
(-1)^{\frac{s}{2}-1}, s\mbox{ even}\\
(-1)^{\frac{s-1}{2}}\cot(\frac{s\pi}{2n})=(-1)^{\frac{s-1}{2}}\gamma_s, s\mbox{ odd}
\end{array}
\right.
\end{array}
\]
The proof will be by induction. For $j=1$ we have $c_1=\xi_1+\cdots+\xi_{n-1}=p_1=\gamma_1$.
Suppose that the proposition is true for all $k<i$.
Then by Newton formula for symmetric functions
\[
je_j=\sum_{k=1}^{j}(-1)^{k-1}p_ke_{j-k}=\sum_{l=1}^{\lceil j/2\rceil}
(p_{2l-1}e_{j-2l+1}-p_{2l}e_{j-2l})
\]
we obtain by writing $c_{j-2l+1}=c_{j-2l}\gamma_{j-2l+1}$
\[
\begin{array}{rl}
je_j=&\displaystyle
\sum_{l=1}^{\lceil j/2\rceil}
\left((-1)^{l-1}\gamma_{2l-1}\gamma_{j-2l+1}-(-1)^{l-1}
\right)c_{j-2l}\\
=&\displaystyle
\sum_{l=1}^{\lceil j/2\rceil}
(-1)^{l-1}(\gamma_{2l-1}\gamma_{j-2l+1}-1)c_{j-2l}\\
\stackrel{*}{=}&\displaystyle
\sum_{l=1}^{\lceil j/2\rceil}
(-1)^{l-1}(\gamma_{2l-1}+\gamma_{j-2l+1})\gamma_jc_{j-2l}\\
=&\displaystyle
\sum_{l=1}^{\lceil j/2\rceil}
(p_{2l-1}c_{j-2l}-p_{2l-2}\gamma_{j-2l+1}c_{j-2l})\gamma_j\ \ (\mbox{here }p_0:=-1)\\
=&\displaystyle
\sum_{l=1}^{\lceil j/2\rceil}
(p_{2l-1}c_{j-2l}-p_{2l-2}c_{j-2l+1})\gamma_j\\
=&\displaystyle
\sum_{l=1}^{\lceil j/2\rceil}
(p_{2l-1}c_{j-1-(2l-1)}-p_{2l-2}c_{j-1-(2l-2)})\gamma_j\\
=&\displaystyle
(-p_0c_{j-1}+\sum_{l=1}^{\lceil (j-1)/2\rceil}
(p_{2l-1}c_{j-1-(2l-1)}-p_{2l}c_{j-1-2l)})\gamma_j\\
\stackrel{**}{=}&\displaystyle
(c_{j-1}+(j-1)c_{j-1})\gamma_j\\
=&\displaystyle
jc_{j-1}\gamma_j=jc_j
\end{array}
\]
Here in $(*)$ we have used the cotangent addition formula
$\cot (\alpha)\cot(\beta)-1=(\cot\alpha+\cot\beta)\cot(\alpha+\beta)$ and in $(**)$
Newton formula for $i-1$ which holds by induction hypothesis. The proposition is thus proved.\qed

For our dihedral configurations we can state the stronger conjecture of Atiyah and Sutcliffe
(\cite{DZb}, Conjecture {\it 2.}) as follows
\begin{equation}
\label{E:Dj1}%                           (*)
n^{\frac{n}{2}}\prod_{k=0}^{n-1}f_k
\geq
2^{\binom{n}{2}}\prod_{i=0}^n(1+\lambda_i^2)^n
\end{equation}
where
\begin{equation}
\label{E:Dj2}%                           (**)
f_k=\sum_{s\geq 0}
\left(
\prod_{j=1}^s\lambda^n_{N-jn-k}\widetilde{E}_{k+sn},\ \ (0\leq k<n)
\right)
\end{equation}
From the factorization
\[
f(y)=h(y)\prod_{i=1}^m(y+\lambda_i)
\]
we can write
\[
\widetilde{E}_k=\sum_{i=0}^{n-1}c_iE_{k-i}
\]
in terms of elementary symmetric functions $E_k=e_k(\lambda_1,\ldots,\lambda_m)$ of our
positive quantities $0<\lambda_1<\cdots<\lambda_m$ with coefficients $c_i$ given in
Proposition \ref{P:2}
(note that $c_0=1\leq c_1\leq\cdots\leq c_{\lfloor\frac{n-1}{2}\rfloor}\geq\cdots\geq c_{n-1}=1$
(unimodality) and $c_i=c_{n-1-i}$ (symmetry)).

Now we shall prove a generalization of the {\DJ}okovi\' c's conjecture which apparently strengthens
(\ref{E:Dj1}).
\begin{theorem}
We have:
\begin{enumerate}
  \item $\displaystyle\prod_{k=0}^{n-1}f_k\geq \prod_{k=0}^{n-1}c_k
  \left(
  \sum_{l=0}^m\left(\prod_{j=0}^{l-1}\lambda_{m-j}E_l\right)
  \right)^n$

  \item$\displaystyle\prod_{k=0}^{n-1}f_k\geq \prod_{k=0}^{n-1}c_k
  \prod_{i=1}^{m}(1+\lambda_{i}^2)^n$
\end{enumerate}
\end{theorem}
\proof
Let us write
\[
f_k=\sum_{l=0}^m\varphi_{kl}E_l
\]
Let us substitute
$\displaystyle\widetilde{E}_{k+sn}=\sum_{i=0}^{n-1}c_iE_{k-i+sn}$
into (\ref{E:Dj2}). Then for fixed $k$ $(0\leq k<n-1)$ and given $l$ $(0\leq l\leq m)$
we seek $s\geq 0$ and $i$, $0\leq i<n$ such that $l=k-1+sn$, i.e.\ $l-k=sn-i$, $0\leq i<n$.
We conclude that $s$ and $i$ are uniquely determined by a division algorithm
(with nonpositive remainder):
\[
s_k:=\left\lceil\frac{l-k}{n}\right\rceil,\ \ \ i_k=s_kn-{l-k}.
\]
Hence
\[
\varphi_{kl}=\prod_{j=1}^{s_k}\lambda_{N-jn-k}^nc_{i_k}
\]
with $s_k$ and $i_k$ just defined. It is easy to see that
\[
s_k=s_0\left(=\left\lceil\frac{l}{n}\right\rceil\right)
\mbox{ and }
i_k=i_0+k\ \mbox{for } 0\leq k\leq n-i_0-1
\]
and
\[
s_k=s_0-1
\mbox{ and }
i_k=i_0+k-n\ \mbox{for } n-i_0\leq k\leq n-1.
\]
\begin{lemma}
For each $l$, $0\leq l\leq m$, we have
\[
\prod_{k=0}^{n-1}\varphi_{kl}=\prod_{j=0}^{l-1}\lambda_{m-j}^n\prod_{j=0}^{n-1}c_j
\]
\end{lemma}
\proof[(of Lemma)]
\[
\begin{array}{l@{\ =\ }l}
\displaystyle
\prod_{k=0}^{n-1}\varphi_{kl}
&\displaystyle
\prod_{k=0}^{n-i_0-1}
\left(\prod_{j=1}^{s_0}\lambda_{N-jn-k}^n\prod_{k=i_0}^{n-1}c_k\right)
\prod_{k=n-i_0}^{n-1}\prod_{j=1}^{s_0-1}\lambda_{N-jn-k}^n\prod_{k=0}^{i_0-1}c_k\\[5mm]
&\displaystyle
\prod_{k=0}^{n-1}\prod_{j=1}^{s_0-1}\lambda_{N-jn-k}^n
\prod_{k=0}^{n-i_0-1}\lambda_{N-s_0n-k}^n\prod_{k=0}^{n-1}c_k
\end{array}
\]
We put now $N=n+m$
\[
\begin{array}{l@{\ =\ }l}
\displaystyle
\hphantom{\prod_{k=0}^{n-1}\varphi_{kl}}
&\displaystyle
\lambda^n_{m}\lambda^n_{m-1}\cdots\lambda^n_{m+n-s_0n-(n-i_0-1)}\prod_{k=0}^{n-1}c_k\\
&\displaystyle
\lambda^n_{m}\lambda^n_{m-1}\cdots\lambda^n_{m-l+1}\prod_{k=0}^{n-1}c_k\\
\end{array}
\]\qed\\
\proof[(of Theorem)]\\
We shall use the H\" older inequality
\[
\begin{array}{rl}
\displaystyle
\prod_{k=0}^{n-1}f_k=\prod_{k=0}^{n-1}
\left(
\sum_{l=0}^{m}\varphi_{kl}E_l
\right)\geq
&\displaystyle
\left(
\sum^{m}_{l=0}
\left(
\prod^{n-1}_{k=1}\varphi_{kl}E_{l}
\right)^{\frac{1}{n}}
\right)^n\\[5mm]
=&\displaystyle
\left(
\sum^{m}_{l=0}
\prod^{l-1}_{j=0}\lambda_{m-j}
\left(
\prod_{j=0}^{n-1}c_j
\right)^{\frac{1}{n}}
E_l
\right)^n\ \ (\mbox{by lemma})\\[7mm]
=&\displaystyle
\left(
\prod_{j=0}^{n-1}c_j\right)
\left(
\sum_{l=0}^{m}
\prod_{j=0}^{l-1}\lambda_{m-j}E_l
\right)^n
\end{array}
\]
Thus {\it 1.} is proved. To obtain {\it 2.} we apply {\DJ}okovi\' c proof of Atiyah conjecture
for type A configurations
\[
\sum_{l=0}^m\prod_{j=0}^{l-1}\lambda_{m-j}E_l\geq\prod_{i=1}^m\left(1+\lambda_i^2\right)
\]
(c.f.\ section 3.)\qed

%%%%%%%%%%%%%%%%%%%%%%%%%%%%%%%%%%%%%%%%%%%%%%%%%%%%%%%%%%%%%%%%%%%%%%%%%%%%%%%%%%%%%%%%%%%%%%%%%%%%%%%%%%

%%%%%%%%%%%%%%%%%%%%%%%%%%%%%%%%%%%%%%%%%%%%%%%%%%%%%%%%%%%%%%%%%%%%%%%%%%%%%%%%%%%%%%
%
% 20060104
%
%%%%%%%%%%%%%%%%%%%%%%%%%%%%%%%%%%%%%%%%%%%%%%%%%%%%%%%%%%%%%%%%%%%%%%%%%%%%%%%%%%%%%%
\section{Appendix}

After the first version of this paper was finished, in the meantime,
we have discovered formulas for the partial derivatives, of the quantities
$\dogstona{n}/\dogstonwoka{n}$,
with respect to variables $\xi_{r}$
(Note that in Theorem \ref{T:1} we have given formulas w.r.t.\ variables $X_r$!).

%%%%%%%%%%%%%%%%%%%%%%%%%%%%%%%%%%%%%%%%%%%%%%%%%%%%%%%%%%%%%%%%%%%%%%%%%%%%%%%%%%%%%%%%%%%%
\begin{lemma}
\label{L:parxi} For $2\leq r\leq n$ the partial derivative w.r.t.\
$\xi_r$ of the quotient $\dogstona{n}/\dogsa{2\ldots n}$ is given by
\[
\left(\dogsa{2\ldots n}\right)^2
\partial_{\xi_r}
\left(\frac{\dogstona{n}}{\dogsa{2\ldots n}}\right)
=
\sum_{i\geq j}
s'_{ij}X_1(X_2\cdots X_j)^2X_{j+1}\cdots X_{i+1}(X_{j+1}-X_{i+2})
\]
where $s'_{ij}$ is the conjugated Schur function
$s_{ij}=s_{ij}(\xi_2,\ldots,\xi_{r-1},\xi_{r+1},\ldots,\xi_n)$ corresponding to a two--rowed partition
$\lambda=(i\geq j)$.\\
In particular for $X_1\geq\cdots\geq X_n>0$ the function $\dogstona{n}/\dogsa{2\ldots n}$
is monotonically increasing w.r.t.\ the variable $\xi_r$ (for $r=1$, too).
\end{lemma}
\proof
By using the formula
$\dogstona{n}=\dogs{1\ldots n-1}{1\ldots\widehat{r}\ldots n}+X_1\xi_r\dogs{2\ldots n}{1\ldots\widehat{r}\ldots n}$
we get
\[
\begin{array}{l}
 \partial_{\xi_r}(\dogstona{n})\dogsa{2\ldots n}- \dogstona{n}\partial_{\xi_r}(\dogsa{2\ldots n})=\\[3mm]
=X_1\dogs{2\ldots n}{1\ldots\widehat{r}\ldots n}
 \left(
 \dogs{2\ldots n-1}{2\ldots\widehat{r}\ldots n}
 +
 X_2\xi_r\dogs{3\ldots n}{2\ldots\widehat{r}\ldots n}
 \right)
 -
 \left(
 \peskd{1}{n-1}{1}{r}{n}+X_1\xi_r\peskd{2}{n}{1}{r}{n}
 \right)X_2\peskd{3}{n}{2}{r}{n}\\[3mm]
=X_1\dogs{2\ldots n}{1\ldots\widehat{r}\ldots n}\dogs{2\ldots n-1}{2\ldots\widehat{r}\ldots n}
 -
 X_2\peskd{1}{n-1}{1}{r}{n}\peskd{3}{n}{2}{r}{n}\\[3mm]
=X_1
\left(
\peskd{2}{n-1}{2}{r}{n}+X_2\xi_1\peskd{3}{n}{2}{r}{n}
\right)
\peskd{2}{n-1}{2}{r}{n}
-
X_2
\left(
\peskd{1}{n-2}{2}{r}{n}+X_1\xi_1\peskd{2}{n-1}{2}{r}{n}
\right)
\peskd{3}{n}{2}{r}{n}\\[3mm]
=X_1
\left(
\peskd{2}{n-2}{2}{r}{n}
\right)^2
-
X_2\peskd{1}{n-2}{2}{r}{n}\peskd{3}{n}{2}{r}{n}
\end{array}
\]
With $e_i=e_i^{(1r)}=e_i(\xi_2,\ldots,\xi_{r-1},\xi_{r+1},\ldots,\xi_n)$ denoting the $i$--th elementary
symmetric function of the truncated alphabet $A^{(1r)}=\{\xi_2,\ldots,\xi_{r-1},\xi_{r+1},\ldots,\xi_n\}$
we have further
\[
\begin{array}{l}
\displaystyle
=X_1
\left(
\sum_{i,j}e_ie_jX_{2\ldots i+1}X_{2\ldots j+1}
\right)
-
X_2
\left(
\sum_{i,j}e_ie_jX_{1\ldots i}X_{3\ldots j+2}
\right)\\[4mm]
\displaystyle
=
\sum_{i,j}e_ie_jX_{1..i+1}X_{2..j+1}
-
\sum_{i,j}e_ie_jX_{1..i}X_{2..j+2}\\[3mm]
\displaystyle
=
\sum_{i,j}
\left|
\begin{array}{cc}
e_i & e_{i+1}\\
e_{j-1} & e_j
\end{array}
\right|
X_{1..i+1}X_{2..j+1}\\[3mm]
\displaystyle
=
\sum_{i\geq j}
\left|
\begin{array}{cc}
e_i & e_{i+1}\\
e_{j-1} & e_j
\end{array}
\right|
X_1(X_{2..j})^2X_{j+1}\cdots X_{i+1}(X_{j+1}-X_{i+2})
\end{array}
\]
Now by Jacobi--Trudy formula we can write
$
\left|
\begin{array}{cc}
e_i & e_{i+1}\\
e_{j-1} & e_j
\end{array}
\right|
$
as the conjugated Schur function $s_{ij}'={s'}_{ij}^{(1r)}$
corresponding to a partition $(i\geq j)$.\qed
%%%%%%%%%%%%%%%%%%%%%%%%%%%%%%%%%%%%%%%%%%%%%%%%%%%%%%%%%%%%%%%%%%%%%%%%%%%%%%%%%%%%%%%%%%%%
\begin{corollary} ($\xi_n$--monotonicity)\\
\label{C:parxi}
We have the following inequality:
\[
\frac{\pes{1}{n}{1}{n}}{\pes{2}{n}{2}{n}}
\geq
\frac{\pes{1}{n-1}{1}{n-1}}{\pes{2}{n-1}{2}{n-1}}
\]
\end{corollary}
\proof
By Lemma \ref{L:parxi} by letting $\xi_n\downarrow 0$ we get
\[
\pes{1}{n}{1}{n}/\pes{2}{n}{2}{n}
\geq \left.\pes{1}{n}{1}{n}/\pes{2}{n}{2}{n}\right|_{\xi_n=0}
=\pes{1}{n-1}{1}{n-1}/\pes{2}{n-1}{2}{n-1}
\]
\qed\\
By using this Corollary we state a strengthening of our Conjecture \ref{C:1}:
\begin{conjecture}
\label{C:1'}%\tag{\ref{C:1}'}
\[
\left(
\pes{1}{n}{1}{n}
\right)^{n-2}
\geq
\pes{2}{n-1}{2}{n-2}
\prod_{k=2}^{n-1}\pesk{1}{k}{n}{1}{k}{n}
\]
\end{conjecture}

%%%%%%%%%%%%%%%%%%%%%%%%%%%%%%%%%%%%%%%%%%%%%%%%%%%%%%%%%%%%%%%%%%%%%%%%%%%%%%%%%%%%%%
%
% 20060105
%
%%%%%%%%%%%%%%%%%%%%%%%%%%%%%%%%%%%%%%%%%%%%%%%%%%%%%%%%%%%%%%%%%%%%%%%%%%%%%%%%%%%%%%
%\begin{remark}
We also have formulas for partial derivative of the quotient
$\dogstona{n}/\dogstonwoka{n}$ w.r.t.\
variable $\xi_r$, $2\leq r\leq n$, which are more complicated than for $k=1$
(given in Lemma \ref{L:parxi}). Without loss of generality we take $r=n$ and
proceed as follows:
\[
\begin{array}{l}
\partial_{\xi_n}(\pes{1}{n}{1}{n})\pesk{1}{k}{n}{1}{k}{n}
-
\pes{1}{n}{1}{n}\partial_{\xi_n}(\pesk{1}{k}{n}{1}{k}{n})
=\\[3mm]
=
X_1\pes{2}{n}{1}{n-1}\pesk{1}{k}{n}{1}{k}{n}
-
X_1\pes{1}{n}{1}{n}\pesk{2}{k}{n}{1}{k}{n-1}\\[3mm]
=
X_1\pes{2}{n}{1}{n-1}
\left(
\pesk{1}{k}{n}{1}{k}{n}+X_1\xi_n\pesk{2}{k}{n}{1}{k}{n-1}
\right)
-
X_1
\left(
\pes{1}{n-1}{1}{n-1}+X_1\xi_n\pes{2}{n}{1}{n-1}
\right)
\pesk{2}{k}{n}{1}{k}{n-1}\\[3mm]
=X_1
\left(
\pes{2}{n}{1}{n-1}\pesk{1}{k}{n-1}{1}{k}{n-1}
-
\pes{1}{n-1}{1}{n-1}\pesk{2}{k}{n}{1}{k}{n-1}
\right)\\[3mm]
=X_1
\left[
\left(
\peskd{2}{n-1}{1}{k}{n-1}+X_2\xi_k\peskd{3}{n}{1}{k}{n-1}
\right)
\pesk{1}{k}{n-1}{1}{k}{n-1}
-\right.\\
\left.\hphantom{=X_1}-\left(
\peskd{1}{n-2}{1}{k}{n-1}+X_1\xi_k\peskd{2}{n-1}{1}{k}{n-1}
\right)
\pesk{2}{k}{n}{1}{k}{n-1}
\right]\\[3mm]
=
X_1
\left[
\peskd{2}{n-1}{1}{k}{n-1}\pesk{1}{k}{n-1}{1}{k}{n-1}
-
\peskd{1}{n-2}{1}{k}{n-1}\pesk{2}{k}{n}{1}{k}{n-1}
+\right.\\
\hphantom{=X_1}
\left.
+\xi_k
\left(
X_2\peskd{3}{n}{1}{k}{n-1}\pesk{1}{k}{n-1}{1}{k}{n-1}
-
X_1\peskd{2}{n-1}{1}{k}{n-1}\pesk{2}{k}{n}{1}{k}{n-1}
\right)
\right]\\[3mm]
=
X_1
\left[
I_1-\xi_kI_2
\right]
\end{array}
\]
Now we first compute
\[
\begin{array}{l}
\displaystyle
I_1=\peskd{2}{n-1}{1}{k}{n-1}\pesk{1}{k}{n-1}{1}{k}{n-1}-\peskd{1}{n-2}{1}{k}{n-1}\pesk{2}{k}{n}{1}{k}{n-1}=\\[3mm]
\displaystyle
\left(
\sum_{i=0}^{k-2}e_iX_{2..i+1}+\sum_{i=k-1}^{n-2}e_iX_{2..i+1}
\right)
\left(
\sum_{j=0}^{k-1}e_jX_{1..j}+\sum_{j=k}^{n-2}e_jX_{1..\widehat{k}..j+1}
\right)-\\[3mm]
\displaystyle
-
\left(
\sum_{j=0}^{k-1}e_jX_{1..j}+\sum_{j=k}^{n-2}e_jX_{1..j}
\right)
\left(
\sum_{i=0}^{k-2}e_iX_{2..i+1}+\sum_{i=k-1}^{n-2}e_iX_{2..\widehat{k}..i+2}
\right)=\\[3mm]
\displaystyle
=
\sum_{i=k-1}^{n-2}\sum_{j=0}^{k-1}e_ie_j\left(X_{2..i+1}X_{1..j}-X_{2..\widehat{k}..i+1}X_{1..j}\right)
+\\[3mm]
\displaystyle
+\sum_{j=k}^{n-2}\sum_{i=0}^{k-2}e_je_i\left(X_{1..\widehat{k}..j+1}X_{2..i+1}-X_{1..j}X_{2..i+1}\right)
+\\[3mm]
\displaystyle
+\sum_{i=k-1}^{n-2}\sum_{j=k}^{n-2}e_ie_j\left(X_{2..i+1}X_{1..\widehat{k}..j+1}-X_{1..j}X_{2..\widehat{k}..i+2}\right)
\end{array}
\]
By replacing, in the middle sum, $j$ with $i+1$ and $i$ with $j-1$, and observing that then
$X_{1..\widehat{k}..i+2}X_{2..j}-X_{1..i+1}X_{2..j}=-(X_{2..i+1}X_{1..j}-X_{2..\widehat{k}..i+2}X_{1..j})$
the contribution of the first two sums is
\[
% ???
%\sum_{0\leq j\leq k-1\leq i\leq n-2}
%???
\sum_{i=k-1}^{n-2}
\sum_{j=0}^{k-1}
\left|
\begin{array}{cc}
e_i & e_{i+1}\\
e_{j-1} & e_j
\end{array}
\right|
X_{2..\widehat{k}..i+1}(X_k-X_{i+2})X_{1..j}
\]
The third sum can similarly be transformed to the following form:
\[
\sum_{k\leq j\leq i\leq n-2}
\left|
\begin{array}{cc}
e_i & e_{i+1}\\
e_{j-1} & e_j
\end{array}
\right|
X_{2..\widehat{k}..i+1}(X_{j+1}-X_{i+2})X_{1..j}
\]
Hence
\[
I_1=\sum_{0\leq j, \max\{j,k-1\}\leq i\leq n-2}
s'_{ij}X_{2..\widehat{k}..i+1}(X_{\max\{j+1,k\}}-X_{i+2})X_{1..j}\ \ (\geq 0)
\]
By a similar manipulation we can obtain the expression for the quantity
\[
\begin{array}{l}
\displaystyle
I_2=X_1\peskd{2}{n-1}{1}{k}{n-1}\pesk{2}{k}{n}{1}{k}{n-1}
-
X_2\peskd{3}{n}{1}{k}{n-1}\pesk{1}{k}{n-1}{1}{k}{n-1}
=\\[3mm]
\displaystyle
=
X_1-X_2
+
\sum_{i=1}^{n-1}\sum_{j\leq min\{k-1,i\}}s'_{ij}X_{2..\widehat{k}..i+2}X_{1..j}(X_{j+1}-X_k)
\geq 0
\end{array}
\]
where $s'_{ij}$ is conjugated Schur function $s'_{ij}={s'}_{ij}^{(kn)}$.
We see that
\[
\left(\dogstonwoka{n}\right)^2\partial_{\xi_n}\left(\frac{\dogstona{n}}{\dogstonwoka{n}}\right)
=
X_1
\left[
I_1-\xi_kI_2
\right]
\]
has both positive and negative terms.
And we have not been able to apply it so far.
%\end{remark}

Now we illustrate use of $\xi$--monotonicity (in addition to $X$--monotonicity)
for proving once more the case $n=4$ of our Conjecture \ref{C:1}:
%By using Corollary \ref{C:2.1} and Corollary \ref{C:parxi} we have:
\[
\begin{array}{ll}
\displaystyle
\frac{(\dogsa{1234})^3}{\dogsa{234}\dogsa{134}\dogsa{124}\dogsa{123}}
=
\frac{\dogsa{1234}}{\dogsa{234}\dogsa{123}}
\frac{\dogsa{1234}}{\dogsa{134}}
\frac{\dogsa{1234}}{\dogsa{124}}
\geq\  (by\ \mbox{$\xi_4$--monotonicity})\\[4mm]
\displaystyle
\geq\frac{1}{\dogsa{23}}
\frac{\dogsa{1234}}{\dogsa{134}}
\frac{\dogsa{1234}}{\dogsa{124}}
\geq\  (by\ \mbox{$X_1$--monotonicity twice and $X_4$--monotonicity})\\[4mm]
\displaystyle
\geq
\frac{1}{\dogsa{23}}
\frac{\dogs{2234}{1243}}{\dogs{234}{143}}
\frac{\dogs{2233}{1234}}{\dogs{223}{124}}
\geq\  (by\ \mbox{$\xi_3$--monotonicity})\\[4mm]
\displaystyle
\geq
\frac{1}{\dogsa{23}}
\frac{\dogs{223}{124}}{\dogs{23}{14}}
\frac{\dogs{2233}{1234}}{\dogs{223}{124}}
=
\frac{\dogs{2233}{1234}}{\dogsa{23}\dogs{23}{14}}
\geq 1
\end{array}
\]

%Since $\dogs{2233}{1234}>\dogsa{23}\dogs{23}{14}$ trivially holds.

Similarly the cases $n=5,6,7$ of Conjecture \ref{C:1} would be, by using
$\xi$--monotonicity and $X$--monotonicity, consequences of the following inequalities
\[
\widetilde{Q}_n\geq 1
\]
where
\[
\begin{array}{l@{\ =\ }l}
\widetilde{Q}_5 &
\dogs{22344}{12345}\dogs{22344}{12345}/
\dogs{234}{234}\dogs{234}{135}\dogs{2244}{1245}\\[3mm]
\widetilde{Q}_6 &
\dogs{223445}{123456}\dogs{233455}{123456}/
\dogs{2345}{2345}\dogs{2345}{1346}\dogs{2345}{1256}\\[3mm]
\widetilde{Q}_7 &
\dogs{2234556}{1234567}\dogs{2334566}{1234567}\dogs{2344566}{1234567}/
\dogs{23456}{23456}\dogs{23456}{13457}\dogs{23456}{12467}\dogs{234566}{123567}
\end{array}
\]
\subsection{Computer verification of the Conjecture \ref{C:1}
(and hence of the Atiyah--Sutcliffe conjecture C3)
for almost collinear $9+1$ configuration.}

Let us now explain our computer verification of the inequality
$\widetilde{Q}_9\geq 1$ where
\[
\widetilde{Q}_9=
\frac{\dogs{223456778}{123456789}\dogs{233456788}{123456789}
\dogs{223456678}{123456789}\dogs{234456788}{123456789}}
{\dogs{2345678}{2345678}\dogs{2345678}{1345679}
\dogs{2345678}{1245689}\dogs{2345678}{1235789}\dogs{22346788}{12346789}}
\]
which refines the case $n=9$ of the Conjecture \ref{C:1}.
We have observed first that $\widetilde{Q}_9$ is symmetric in partial alphabets
\[
A_1=\{\xi_1,\xi_2,\xi_8,\xi_9\},\
A_2=\{\xi_3,\xi_4,\xi_6,\xi_7\},\
A_3=\{\xi_5\}
\]
then by introducing the elementary symmetric functions $\{e_1, e_2, e_3, e_4\}$
of $A_1$ and $\{f_1, f_2, f_3, f_4\}$ of $A_2$ we first computed the products
\[
\dogs{2345678}{2345678}\dogs{2345678}{1345679}
\mbox{ and }
\dogs{2345678}{1245689}\dogs{2345678}{1235789}
\]
in terms of $\{e_1, e_2, e_3, e_4,f_1, f_2, f_3, f_4,\xi_5\}$.
Then by successive application of Stembridge's {\tt Maple} {\tt SF} package we expressed the difference
$\Delta:=numer(\widetilde{Q}_9)-denom(\widetilde{Q}_9)$ in terms of the Schur
functions of both alphabets $A_1$ and $A_2$.
Then we factored each coefficient in such a multi--Schur expansion and into nonmonomial
factors we substituted $X_2=X_3+h_2$, $X_3=X_4+h_3$, $\ldots$, $X_7=X_8+h_7$. Then the
computation showed that the coefficients of all monomials in $X_8,h_2,\ldots,h_7$
were nonnegative. The factoring out the trivial monomial factors in $X_2,\ldots,X_8$
(which are trivially nonnegative) was crucial because otherwise the expansion
of multi--Schur function coefficients in terms of increments $h_2,\ldots,h_7$ may
not be feasible.

%%%%%%%%%%%%%%%%%%%%%%%%%%%%%%%%%%%%%%%%%%%%%%%%%%%%%%%%%%%%%%%%%%%%%%%%%%%%%%%%%%%%%%%%%%%%%%%%%%%%%%%%%%

\section{Appendix 2}

Here we first recall a remarkable inequality of I.\ Schur (c.f.\
J.~Michael Steele: The Cauchy--Schwarz Master Class, Cambridge
University Press, 2004.)

For all values $x,y,z\geq 0$ and all $\alpha\geq 0$ we have
\[
\begin{array}{l}
I_\alpha(x,y,z):=
\sum x^\alpha(x-y)(x-z)=\\
\hphantom{I_\alpha(x,y,z)\ }=
x^\alpha(x-y)(x-z)+y^\alpha(y-x)(y-z)+z^\alpha(z-x)(z-y)\geq 0
\end{array}
\]
with equality iff either $x=y=z$ or two of the variables are equal and the third is zero.
Note that $I_\alpha$ is a symmetric function. For a proof we can assume $0\leq x\leq y\leq z$.
Then clearly $x^\alpha(x-y)(x-z)\geq 0$ and by grouping the other two terms we get
$(z-y)[z^\alpha(z-x)-y^\alpha(y-x)]\geq 0$ by observing that $z\geq y$ and $z-x\geq y-x$.

Now we state and prove several properties of a function
\[
d_3(x,y,z):=(x+y-z)(x-y+z)(-x+y+z)\ \ (x,y,z\geq 0)
\]
which frequently appears in the main part of this paper.

We note that the area $A=A(a,b,c)$ of a triangle with sides lengths $a,b,c$
is given, according to the Heron-s formula:
\[
\begin{array}{l@{\ =\ }l}
(4A)^2
& (a+b+c)(a+b-c)(a-b+c)(-a+b+c)\\[2mm]
& (a+b+c)d_3(a,b,c)\\[2mm]
& 2a^2b^2+2a^2c^2+2b^2c^2-a^4-b^4-c^4
\end{array}
\]

Properties of the function $d_3$:\\
\begin{proposition}
\label{P:A2}
We have the following identities and inequalities:
\begin{enumerate}
  \item $xyz-d_3(x,y,z)=\sum x(x-y)(x-z)\geq 0$

  \item $d_3(x,y,z)^2-d_3(x^2,y^2,z^2)=$\\
  $=\sum x^2(y^2-yz+z^2-x^2)^2+\left(\sum x(y^2-yz+z^2-x^2)\right)^2\geq 0$

  \item $d_3(x,y,z)^2-d_3(x^2,y^2,z^2)=$\\
  $=8x^2y^2z^2-2(xyz+x^3+y^3+z^3)d_3(x,y,z)\geq 0$

  \item $(x+y+z)^2d_3(x,y,z)^2-3(x^2+y^2+z^2)d_3(x^2,y^2,z^2)=$\\
  $=4\sum x^4(x^2-y^2)(x^2-z^2)\geq 0$

  \item $(x+y+z)(X+Y+Z)d_3(x,y,z)d_3(X,Y,Z)-3(xX+yY+zZ)d_3(xX,yY,zZ)=$\\
  $2\sum (x^2(x^2-y^2)X^2(X^2-Z^2)+X^2(X^2-Y^2)x^2(x^2-z^2))+
   (x^2(Y^2-Z^2)+y^2(Z^2-X^2)+z^2(X^2-Y^2))^2\geq 0$
\end{enumerate}
\end{proposition}
\proof

All identities {\it 1.}--{\it 5.\ } can be easily checked by expansion. The inequality in
{\it 1.\ } follows from Schur's inequality ($\alpha=1$), in {\it 2.\ } it is evident since
the rhs is the sum of four squares (see \cite{Carroll}).
Case {\it 3.\ } follows from {\it 2.\ } Case {\it 4.\ }
follows from Schur's inequality ($\alpha=2$). Case {\it 5.\ } follows from a generalization
of the case $\alpha=2$ of Schur's inequality:
\[
\begin{array}{l}
I\! I_2(x,y,z,X,Y,Z)=\sum x(x-y)X(X-Z)=\\
= x(x-y)X(X-Z)+y(y-x)Y(Y-Z)+z(z-x)Z(Z-Y)\geq 0
\end{array}
\]
(by letting $y=x+h$, $z=y+k$, $Y=X+H$, $Z=Y+K$,\ $h,k,H,K\geq 0$).\qed

\begin{corollary}
From the Proposition we get the following inequalities:
\[
d_3(x,y,z)\leq xyz\ \ (from\ \mbox{{\it 1.\ }})
\]
and a stronger inequality $d_3(x,y,z)\leq 4x^2y^2z^2/(xyz+x^3+y^3+z^3)$ (from {\it 3.})

 From {\it 2.\ } we have the inequality
\[
d_3(x,y,z)^2\geq d_3(x^2,y^2,z^2)
\]
which can also be obtained from {\it 4.\ } (which implies famous Finsler--Hadwiger inequality)
by using the inequality $(x+y+z)^2\leq 3(x^2+y^2+z^2)$.

The inequality {\it 5.\ }, with the help of Chebyshev inequality
\[
(x+y+z)(X+Y+Z)\leq 3(xX+yY+zZ)\ \ (x\leq y\leq z,\ \ X\leq Y\leq Z)
\]
gives us the following inequality (which seems to be new):
\[
d_3(x,y,z)d_3(X,Y,Z)\geq d_3(xX,yY,zZ)
\]
(when $0\leq x\leq y\leq z$, $0\leq X\leq Y\leq Z$).
\end{corollary}
\begin{remark}
If $a, b, c$ are side lengths of a triangle then inequality
$d_3(a,b,c)\leq abc$ follows also directly from the following identity
\[
abc-d_3(a,b,c)=\frac{1}{2}[(-a+b+c)(b-c)^2+(a-b+c)(a-c)^2+(a+b-c)(a-b)^2]
\]
from which we also have the following inequality
\[
8(abc-d_3(a,b,c))^3\geq d_3(a,b,c)(a-b)^2(a-c)^2(b-c)^2
\]
\end{remark}

%%%%%%%%%%%%%%%%%%%%%%%%%%%%%%%%%%%%%%%%%%%%%%%%%%%%%%%%%%%%%%%%%%%%%%%%%%%%%%%%%%%%%%%%%%%%%%%%%%%%%%%%%%

%%%%%%%%%%%%%%%%%%%%%%%%%%%%%%%%%%%%%%%%%%%%%%%%%%%%%%%%%%%%%%%%%%%%%%%%%%%%%%%%%%%%%%%%%%%%%%%%%%%%%%%%%%


\begin{thebibliography}{9}

\bibitem{MA1}
M. Atiyah.
\newblock The geometry of classical particles.
\newblock {\em Surveys in Differential Geometry}
(International Press) {\bf 7} (2001).

\bibitem{MA2}
M. Atiyah.
\newblock Configurations of points.
\newblock {\em Phil. Trans. R. Soc. Lond. A}
{\bf 359} (2001), 1375--1387.

\bibitem{AS1}
M. Atiyah and P. Sutcliffe.
\newblock Polyhedra in Physics, Chemistry and Geometry,
\newblock {\em Milan Journal of Mathematics}, Vol 71, Number 1/September 2003, 33--58

\bibitem{AS}
M. Atiyah and P. Sutcliffe.
\newblock The geometry of point particles.
\newblock {\em Royal Society of London Proceedings Series A},
\newblock vol. 458, Issue 2021., 1089.--1116.

\bibitem{Carroll}
C.\ E.\ Carroll, C.\ C.\ Yang and S.\ Ahn,
\newblock Some triangle inequalities and generalizations,
\newblock {\em Canad.\ Math.\ Bull.\ Vol.\ 23(3)}, 1980, 267--274

\bibitem{EN}
M. Eastwood and P. Norbury,
\newblock A proof of Atiyah's conjecture on configurations
of four points in Euclidean three-space.
\newblock {\em Geometry \& Topology} {\bf 5} (2001), 885--893.

\bibitem{DZ}
D.\v{Z}. {\DJ}okovi\'{c}, Proof of Atiyah's conjecture for two special
types of configurations, arXiv:math.GT/0205221 v4, 11 June 2002.
{\em Electron. J. Linear Algebra} {\bf 9} (2002), 132--137.

\bibitem{DZb}
{\DJ}okovi\'{c}, D. \v{Z}., Verification of Atiyah's conjecture for some
nonplanar configurations with dihedral symmetry.In \emph{ Publ.
Inst. Math., Nouv. Sér.} 72(86),  (2002) 23--28.

%\bibitem{DZb}
%D.\v{Z}. {\DJ}okovi\'{c}, Verification of Atiyah's conjecture for
%some nonplanar configurations with dihedral symmetry,
%arXiv:math.GT/0208089 v2, 13 Aug 2002.

\bibitem{Jacobi}
I.\ G.\ Macdonald
\newblock {\it Symmetric functions and Hall polynomials}
\newblock 2$^{\mbox{nd}}$ edition,
\newblock Oxford University Press,
\newblock 1995.

\bibitem{Moebius}
A.\ F.\ M\" obius
\newblock \" Uber eine Methode, um von Relationen, welche der Longimetrie
angeh\" oren, zu entsprechenden S\" atzen der Planimetrie zu
gelangen,
\newblock Ber.\ Verhandl.\ k\" onigl.\ S\" achs.\ Ges.\ Wiss.\
Leipzig Math.\ Phys.\ Cl.,
\newblock 1852., 41--45.

\bibitem{DSIU}
D.\ Svrtan, I.\ Urbiha
\newblock Atiyah-Sutcliffe Conjectures for Almost Collinear Configurations and Some New Conjectures for Symmetric Functions,
\newblock math.AG/0406386

\bibitem{Wolstenholme}
J. Wolstenholme
, A Book of Mathematical Problems on Subjects
Included in the Cambridge Course, London and Cambridge, 1867, 56

\bibitem{Walker}
A.W.Walker, 'Problem 300', Nieuw Arch. Wisk. (3) 19 (1971), 224.

\end{thebibliography}
\end{document}